\documentclass[brochure,12pt,english]{bourbaki}
\usepackage[utf8]{inputenc}
\usepackage{eucal}
\usepackage[obeyspaces]{url}
\usepackage{epsfig}

\theoremstyle{plain}
\newtheorem{quest}[defi]{Question}
\newtheorem{hypot}[defi]{Assumption}

\newcommand{\tmu}{\tilde{\mu}}
\newcommand{\sop}{\Sigma}
\newcommand{\T}[1]{{}^t{{#1}}}

\newcommand{\Cc}{\mathbf{C}}

\newcommand{\Zz}{\mathbf{Z}}

\newcommand{\Rr}{\mathbf{R}}

\newcommand{\Qq}{\mathbf{Q}}
\newcommand{\Fp}{\mathbf{F}}

\newcommand{\orbit}{\mathcal{O}}
\newcommand{\satur}{r}
\newcommand{\sifted}{\mathcal{S}}
\DeclareMathOperator{\SL}{SL}
\DeclareMathOperator{\GL}{GL}
\DeclareMathOperator{\Sp}{Sp}

\DeclareMathOperator{\SO}{SO}

\newcommand{\apol}{\mathcal{A}}
\newcommand{\curv}{\mathcal{C}}
\newcommand{\lubotz}{L}

\def\stacksum#1#2{{\stackrel{{\scriptstyle #1}}
{{\scriptstyle #2}}}}
\def\multsum{\mathop{\sum\cdots\sum}\limits}

\newcommand{\mods}[1]{\ (\mathrm{mod}\ {{#1}})}
\newcommand{\uple}[1]{\text{\boldmath${#1}$}}
\newcommand{\eps}{\varepsilon}
\renewcommand{\leq}{\leqslant}
\renewcommand{\geq}{\geqslant}

\newcommand{\ra}{\rightarrow}
\newcommand{\lra}{\longrightarrow}

\newcommand{\fleche}[1]{\stackrel{#1}{\lra}}

\newcommand{\proba}{\mathbf{P}}

\date{November 2010} \bbkannee{63\`eme ann\'ee, 2010-2011}
\bbknumero{1028} \title{Sieve in expansion}
\author{Emmanuel KOWALSKI}
\address{ETH Z\"urich -- DMATH\\
R\"amistrasse 101\\
8092 Z\"urich, Switzerland}
\email{kowalski@math.ethz.ch}
\begin{document}
\maketitle


\section{Introduction}\label{sec-intro}

This report presents recent works extending sieve methods, from their
classical setting, to new situations characterized by the targeting of
sets with exponential growth, arising often from discrete groups like
$\SL_m(\Zz)$ or sufficiently big subgroups. 
\par
A recent lecture of Sarnak~\cite{sarnak-lecture} mentions some of the
original motivation (related to the Markov equation and closed
geodesics on the modular surface). The first general results
concerning these sieve problems appeared around $2005$ in preprint
form, and Bourgain, Gamburd and Sarnak have written a basic paper
presenting its particular features~\cite{bgs}.  Other applications,
with a very different geometric flavor, also appeared independently
around that time, first (somewhat implicitly) in some works of
Rivin~\cite{rivin}.
\par
The most crucial feature in applying sieve to these new situations is
their dependence on \emph{spectral gaps}, either in a discrete setting
(related to expander graphs or to Property $(\tau)$ of
Lubotzky~\cite{lubotzky}) or in a geometric setting (generalizing for
instance Selberg's result that $\lambda_1\geq 3/16$ for the spectrum
of the Laplace operator on the classical hyperbolic modular surfaces).
\par
The outcome of these developments is that there now exist very general
sieve inequalities involving, roughly speaking, discrete objects with
exponential growth. Moreover, their applicability (including to
problems seemingly unrelated with classical analytic number theory, as
we will show) has expanded enormously, as -- partly motivated by these
new applications of sieve methods -- many new cases of spectral gaps
have become available. Particularly impressive are the results on
expansion in finite linear groups (due to many people, but starting
from the breakthrough of Helfgott~\cite{helfgott} for $\SL_2$), and
those concerning applications of ergodic methods to lattices in
semisimple groups with Property $(\tau)$ (developed most generally by
Gorodnik and Nevo~\cite{gn}).
\par
Before going towards the heart of this report, we state here a
particularly concrete and appealing result arising from sieve in
expansion. We recall first that $\Omega(n)$ is the arithmetic function
giving the number of prime factors, counted with multiplicity, of a
non-zero integer $n$, extended so that $\Omega(0)=+\infty$.

\begin{theo}\label{th-bgs-ex}
  Let $\Lambda\subset \SL_m(\Zz)$ be a Zariski-dense subgroup, for
  instance the group $\lubotz$ generated by the elements
\begin{equation}\label{eq-lub}
\begin{pmatrix}
1& \pm 3\\
0& 1
\end{pmatrix},
\quad
\begin{pmatrix}
1 & 0\\
\pm 3& 1
\end{pmatrix}\in \SL_2(\Zz),
\end{equation}
in the case $m=2$
\par
Let $f$ be an integral polynomial function on $\Zz^m$, which is
non-constant. Let $x_0\in\Zz^m-\{0\}$ be a fixed vector. There
exists an integer $\satur=\satur(f,x_0,\Lambda)\geq 1$ such that
the set
$$
\orbit_f(x_0;\satur)= \{ \gamma\in \Lambda\,\mid\,
\Omega(f(\gamma\cdot x_0))\leq \satur \}
$$
is \emph{Zariski-dense} in $\SL_m$, and in particular is infinite. In
fact, there exists such $\satur$ for which $\orbit_f(x_0;\satur)$ is
\emph{not thin}, in the sense of~\cite[Def. 3.1.1]{serre-galois}.
\end{theo}

Part of the point, and it will be emphasized below, is that $\Lambda$
may have infinite index in $\SL_m(\Zz)$. In particular, for $m=2$,
this is the case for the group $L$ generated by the
matrices~(\ref{eq-lub}).
\par
\medskip
\textbf{Notation.} We recall here some basic notation.
\par
\noindent -- The letter $p$ will always refer to prime numbers; for a
prime $p$, we write $\Fp_p$ for the finite field $\Zz/p\Zz$, and we
write $\Fp_q$ for a field with $q$ elements. For a set $X$, $|X|$ is
its cardinality, a non-negative integer or $+\infty$.
\par
\noindent -- The Landau and Vinogradov notation $f=O(g)$ and $f\ll g$
are synonymous, and $f(x)=O(g(x))$ for all $x\in D$ means that there
exists an ``implied'' constant $C\geq 0$ (which may be a function of
other parameters, explicitly mentioned) such that $|f(x)|\leq Cg(x)$
for all $x\in D$.  This  definition \emph{differs} from that of
N. Bourbaki~\cite[Chap. V]{FVR} since the latter is of topological
nature. On the other hand, the notation $f(x)\sim g(x)$ and $f=o(g)$
are used with the asymptotic meaning of loc. cit.
\par
\medskip
\par
\textbf{Acknowledgments.} Thanks are due to J. Bourgain, N. Dunfield,
E. Fuchs, A. Gamburd, F. Jouve, A. Kontorovich, H. Oh, L. Pyber,
P. Sarnak, D. Zywina and others for their help, remarks and
corrections concerning this report. In particular, discussions with
O. Marfaing during his preparation of a Master Thesis on this
topic~\cite{marfaing} were very helpful.

\section{Motivation}\label{sec-motivation}

Sieve methods are concerned with multiplicative properties of sets of
integers. Thus to expand the range of the sieve, one should describe
new sets of integers to investigate. To present the spirit of this
survey, we first give two examples of such sets of integers, which are
rather unusual from a sieve perspective. One of them is a particularly
appealing instance of ``sieve in orbits'', first considered
in~\cite{bgs}: the distribution of curvatures of integral Apollonian
circle packings. The second is even more surprising to look at: it has
to do with the first homology of certain ``random''
$3$-manifolds. Although we will say rather less about it later on, it
presents some unusual features, and suggests interesting questions.

\subsection{Apollonian circle packings}\label{ssec-apol}

It is a very classical geometrical fact that, given three circles
$(\bigcirc_1,\bigcirc_2,\bigcirc_3)$ in the plane which are pairwise
tangent to each other, and have disjoint ``interiors'' (the discs they
bound), with radii $(r_1,r_2,r_3)$ and curvatures
$(c_1,c_2,c_3)=(r_1^{-1},r_2^{-1},r_3^{-1})$, one can find two more
circles (say $(\bigcirc_4,\bigcirc'_4)$ with curvatures $(c_4,c'_4)$),
so that both
$$
(\bigcirc_1,\bigcirc_2,\bigcirc_3,\bigcirc_4)\text{ and }
(\bigcirc_1,\bigcirc_2,\bigcirc_3,\bigcirc'_4)
$$
are four pairwise tangent circles (with disjoint ``interiors''). In
fact, this applies also to negative radii or curvatures, where a
negative radius is interpreted to mean that the ``interior'' of the
circle should be the complement of the bounded disc (see Figure~1).
Such $4$-tuples are called \emph{Descartes configurations}, since a
result of Descartes states that the two sets of four curvatures
satisfy the quadratic equations
$$
Q(c_1,c_2,c_3,c_4)=Q(c_1,c_2,c_3,c'_4)=0
$$
where
$$
Q(x,y,z,t)=2(x^2+y^2+z^2+t^2)-(x+y+z+t)^2.
$$
\par
In particular, if $(\bigcirc_1,\bigcirc_2,\bigcirc_3,\bigcirc_4)$ are
such that their curvatures are all \emph{integers}, then we obtain an
integral quadratic equation for $c'_4$ where one solution (namely,
$c_4$) is an integer: thus $c'_4$ is an integer also. Moreover, again
given $(\bigcirc_1,\bigcirc_2,\bigcirc_3,\bigcirc_4)$ with integral
curvatures, there are also circles
$$
\bigcirc'_1,\bigcirc'_2,\bigcirc'_3,
$$
for which, for instance, the circles
$$
(\bigcirc'_1,\bigcirc_2,\bigcirc_3,\bigcirc_4)
$$
form a Descartes configuration, and as above, the curvatures $c'_1$,
$c'_2$, $c'_3$ are integers. In fact, these new curvatures are given
-- by solving the quadratic equation using the known root -- as
\begin{align*}
(c'_1,c_2,c_3,c_4)=(c_1,c_2,c_3,c_4)\cdot \T{s_1},\\
(c_1,c'_2,c_3,c_4)=(c_1,c_2,c_3,c_4)\cdot \T{s_2},\\
(c_1,c_2,c'_3,c_4)=(c_1,c_2,c_3,c_4)\cdot \T{s_3},\\
(c_1,c_2,c_3,c'_4)=(c_1,c_2,c_3,c_4)\cdot \T{s_4},
\end{align*}
where the matrices $s_1$, \ldots, $s_4$ are in the group $O(Q,\Zz)$ of
integral automorphisms of the quadratic form above, namely
$$
s_1=\begin{pmatrix}
-1&2&2& 2\\
&1&&\\
&&1&\\
&&&1
\end{pmatrix},\quad\quad
s_2=\begin{pmatrix}
1&&&\\
2&-1& 2&2\\
&&1&\\
&&&1
\end{pmatrix},\quad\quad
$$
and $s_3$ and $s_4$ are similar. Note that $s_i^2=1$ for all $i$, and
one can in fact show that these are the only relations satisfied by
those matrices.
\par
Each of the new sets of curvatures can be used to iterate the process;
in other words, denoting by $\apol$ the subgroup of $O(Q,\Zz)$
generated by the $s_i$, the integers arising as coefficients of a
vector in the orbit $\apol\cdot \uple{c}$ of a ``root quadruple''
$\uple{c}=(c_1,\ldots, c_4)$, are all the curvatures of circles
arising in this iterative circle packing. These are \emph{Apollonian
  circle packings}, and the first step is described in Figure~1 in one
particular case (where one notices the convention dealing with
negative curvatures).


\begin{figure}\label{fig-pk}
\begin{center}
\epsfig{file=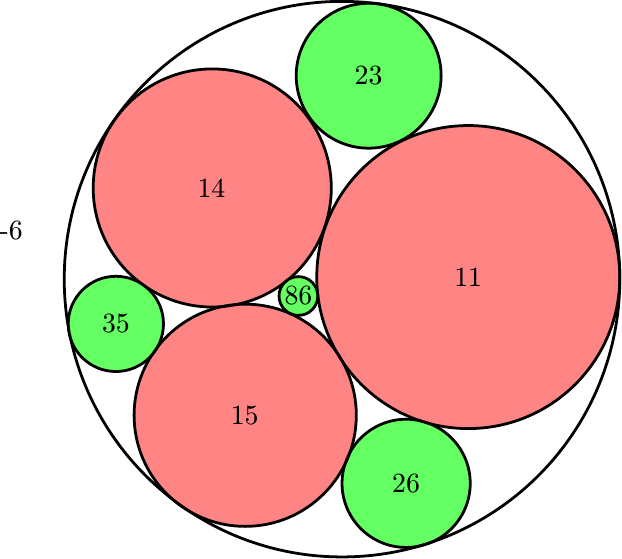,width=0.5\textwidth}
\end{center}
\caption{Apollonian circle packing for $\uple{c}=(-6,11,14,15)$; the
  labels are the curvatures.}
\end{figure}
\par
The set $\curv(\uple{c})$ of these curvatures, considered with or
without multiplicity, is our first example of integers to sieve
for. It is clear from the outset that such an attempt will be deeply
connected with the understanding of the group $\apol$. Moreover, an
interest in the multiplicative properties of the elements of
$\curv(\uple{c})$ will obviously depend on the properties of the
reduction maps
$$
\apol\ra \apol_p=\apol \mods{p}\subset O(Q,\Zz/p\Zz),
$$
modulo primes, and in particular in the image of this reduction map.
\par
The following features of $\apol$ illustrate a basic property that
makes the question challenging:
\par
-- The group $\apol$ is ``big'' in some sense: it is Zariski-dense in
$O(Q)$ (as a $\Qq$-algebraic group), so that, if one can only use
polynomial constructions, $\apol$ is indistinguishable from the big
Lie group $O(Q,\Cc)$.
\par
-- However, $\apol$ is ``small'' in some other sense: specifically,
$\apol$ has \emph{infinite index} in $O(Q,\Zz)$. Stated differently,
the quotient $\apol\backslash O(Q,\Rr)$ (a three-manifold) has
infinite volume for its natural measure, induced from Haar measure on
$O(Q,\Rr)$.

\begin{rema}
Arithmetic properties of Apollonian packings were first
discussed in~\cite{glmwy}, including some properties of
$\curv(\uple{c})$; the use of sieve to study $\curv(\uple{c})$ was
begun in~\cite{bgs}.
\end{rema}

\subsection{Dunfield-Thurston random $3$-manifolds}
\label{ssec-geom}

Our second example of integers to sieve from is chosen partly as an
illustration of the great versatility of sieve, and partly because it
will lead to some interesting examples later. It is based on a paper
of Dunfield and Thurston~\cite{dt} (which did not explicitly introduce
sieve); related work is due to Maher~\cite{maher} and the
author~\cite{cup} (where sieve is explicitly present).
\par
Let $g\geq 2$ be a fixed integer. Fix also a handlebody $H_g$ of genus
$g$; it is a connected oriented compact $3$-manifold with boundary,
and this boundary surface $\Sigma_g$ is a surface of genus $g$
(compact connected and oriented); in other words, for $g=2$, $H_g$ is
a solid double ``doughnut'', and $\Sigma_g$ its $2$-dimensional
boundary. A very classical way (going back to Heegaard) of
constructing compact $3$-manifolds (also connected, oriented, without
boundary) is the following: take a homeomorphism $\phi$ of $\Sigma_g$,
and consider the manifold
$$
M_{\phi}=H_g \cup_{\phi} H_g
$$
obtained by gluing two copies of $H_g$ using the map $\phi$ to
identify points on their common boundary.
\par
As may seem intuitively reasonable, the manifold $M_{\phi}$ does not
change if $\phi$ is changed continuously; this means that $M_{\phi}$
really depends only on the mapping class of $\phi$ in the mapping
class group $\Gamma_g$ of $\Sigma_g$ (roughly, the group of
``discrete'' invariants of the homeomorphisms of the surface; these
groups, as was recently discussed in this seminar~\cite{paulin}, have
many properties in common with the groups $\SL_m(\Zz)$, or better with
their quotients $\Sp_{2g}(\Zz)$).
\par
The topic of~\cite{dt} (which is partly inspired by the Cohen-Lenstra
heuristic for ideal class groups of number fields) is the
investigation of the statistic properties of the fundamental group
$\pi_1(M_{\phi})$ when $\phi$ is taken as a ``random'' elements of
$\Gamma_g$ (in a sense to be described precisely below), in particular
the study of the abelianization $H_1(M_{\phi},\Zz)$ of
$\pi_1(M_{\phi})$, motivated by the virtual Haken Conjecture
(according to which every compact $3$-manifold with infinite
fundamental group should have a finite covering $N$ such that the
abelianization of $\pi_1(N)$ is infinite).
\par
Motivated by this, our second example of set of integers is, roughly,
the set of the orders of torsion subgroups of $H_1(M_{\phi},\Zz)$, as
$\phi$ runs over $\Gamma_g$. Or rather, since here multiplicity is
very hard to control, one should think of this as the map
$$
\Gamma_g\ra |H_1(M_{\phi},\Zz)| \in \{0,1,2,3,\ldots \}\cup 
\{+\infty\}.
$$
\par
It may be fruitful to think about these integers from a sieve point of
view because of the ``local'' information given by the homology with
coefficients in $\Fp_p$ for $p$ prime, namely
$$
H_1(M_{\phi},\Zz)\otimes \Zz/p\Zz=H_1(M_{\phi},\Zz/p\Zz),
$$
and the sieve-like description
$$
\dim H_1(M_{\phi},\Zz)\otimes \Qq\geq 1
\Longleftrightarrow (\text{For all primes $p$, }
\dim_{\Zz/p\Zz} H_1(M_{\phi},\Zz/p\Zz)\geq 1)
$$
of the manifolds with infinite first homology (this is because
$H_1(M_{\phi},\Zz)$ is a finitely-generated group).
\par
A certain similarity with the previous example arises here from the
fact that, as shown by Dunfield and Thurston, there is a natural
isomorphism
\begin{equation}\label{eq-homology}
H_1(M_{\phi},\Zz)\simeq V/\langle J,\phi_* J\rangle
\end{equation}
where $V=H_1(\Sigma_g,\Zz)\simeq \Zz^{2g}$ is the first homology of
the surface $\Sigma_g$, $J$ is the image in $V$ of $H_1(H_g,\Zz)\simeq
\Zz^g$, which is a (fixed!) Lagrangian subspace for the intersection
pairing on $V$, and $\phi_*$ denotes the induced action of $\phi$ on
$V$. Thus $H_1(M_{\phi},\Zz)$ only depends on $\phi_*$, which is an
element of the discrete group $\Sp(V)\simeq \Sp_{2g}(\Zz)$ of
symplectic maps on $V$, with respect to the intersection
pairing. Moreover, the reduction modulo $p$ is given by
\begin{equation}\label{eq-homologyp}
 H_1(M_{\phi},\Zz/p\Zz)\simeq V_p/\langle J_p,\phi_{*} J_p\rangle
\end{equation}
where $V_p=V/pV$, $J_p=J/pJ$, and therefore it depends only on the
reduction modulo $p$ of $\phi_*$, an element of the finite group
$\Sp(V/pV)\simeq \Sp_{2g}(\Zz/p\Zz)$.

\section{A quick survey of sieve}
\label{sec-sieve}

In this section, we survey quickly some of the basic principles of
sieve methods, and state one version of the fundamental result that
evolved from V. Brun's first investigations. Our goal is to make the
sieve literature, and its terminology and notation, accessible to
non-experts. In particular, this section is essentially
self-contained; only in the next one do we start to fit the examples
above (and many others) in the sieve framework.

\subsection{Classical sieve}

The classical sieve methods arose from natural questions related to
the way multiplicative constraints on positive integers (linked to
restrictions on their prime factorization) can interact with additive
properties. The best known among these questions, and a motivating one
for V. Brun and many later arithmeticians, is whether there exist
infinitely many prime numbers $p$ such that $p+2$ is also prime, but
the versatility of sieve methods is quite astounding. For many
illustrations, and for background information, we refer to the very
complete modern treatment found in the recent book of J. Friedlander
and H. Iwaniec~\cite{FI}.
\par
In the usual setting, the basic problem of sieve theory is the
following: given a sequence $\mathcal{F}=(a_n)_{n\geq 1}$ of
non-negative real numbers (usually with a finite support, which is
supposed to be a parameter tending to infinity) and a (fixed) subset
$\mathcal{P}$ of the primes (e.g., all of them), one seeks to
understand the sum
$$
S(\mathcal{F},z)=\sum_{\stacksum{n\geq 1}{(n,P(z))=1}}{a_n},\quad
\text{ where }\quad P(z)=\prod_{\stacksum{p\in\mathcal{P}}{p<z}}{p},
$$
which encodes the contribution to the total sum
$$
S(\mathcal{F})=\sum_{n\geq 1}{a_n},
$$
of the integers not divisible by the primes in $\mathcal{P}$ which are
$<z$, and one wishes to do so using properties of the given sequence
which are encoded in \emph{sieve axioms} or \emph{sieve properties}
concerning the congruence sums
\begin{equation}\label{eq-cong-sums}
S_d(\mathcal{F})=\sum_{\stacksum{n\geq 1}{n\equiv 0\mods{d}}}{a_n}.
\end{equation}
\par
The fundamental relation between these quantities is the well-known
\emph{inclusion-exclusion} formula\footnote{\ Where $\mu(d)$ is the
  M\"obius function, $0$ for non-squarefree integers, and otherwise
  equal to $(-1)^{\Omega(n)}$.}
$$
S(\mathcal{F},z)=\sum_{d\mid P(z)}{\mu(d)S_d(\mathcal{F})},
$$
and the basic philosophy is that, for many sequences of great
arithmetic interest, the congruence sums above can be understood quite
well. Indeed, the notion of ``sieve of dimension $\kappa>0$'' arises
as corresponding to a sequence $\mathcal{F}$ for which, intuitively,
the ``density'' of the sequence over integers divisible by a prime
number $p$ (in $\mathcal{P}$) is approximately $\kappa p^{-1}$, at
least on average over $p$, i.e., we have
\begin{equation}\label{eq-sieve-axiom}
S_d(\mathcal{F})=g(d)S(\mathcal{F})+r_d(\mathcal{F}),
\end{equation}
where $r_d(\mathcal{F})$ is considered as a ``small'' remainder and
$g$ is a multiplicative function of $d\geq 1$ for which $g(p)$
satisfies
\begin{equation}\label{eq-dimension2}
g(p)=\frac{\kappa}{p}+O(p^{-1-\delta})
\end{equation}
for some $\delta>0$ (or even weaker or averaged versions of this, such
as
\begin{equation}\label{eq-dimension}
\sum_{\stacksum{p\leq x}{}}{g(p)\log p}=\kappa\log x+O(1)\ ;
\end{equation}
since we have
$$
\sum_{p\leq x}{\frac{\log p}{p}}=\log x+O(1),
$$
for $x\geq 2$, by the Prime Number Theorem, such an assumption is
consistent with the heuristic suggested above).
\par
This dimension condition often means that the sum $S(\mathcal{F},z)$
corresponds to the number of integers (in a finite sequence) which,
modulo the primes $p\in\mathcal{P}$, must avoid  $\kappa$ residue
classes. 

\begin{exem}\label{ex-1}
  A characteristic example is the sequence
  $\mathcal{F}_f=\mathcal{F}_{f,X}$, associated with a fixed monic
  polynomial $f\in\Zz[T]$ of degree $r\geq 1$ and a (large) parameter
  $X$, defined as the multiplicity
\begin{equation}\label{eq-olda}
a_n=|\{m\leq X\,\mid\, f(m)=n\}|
\end{equation}
of the representations of an integer as a value $f(m)$ with $m\leq
X$. In that case, if $\mathcal{P}$ is the set of all primes and
$z\approx X^{1/2r}$, it follows that
\begin{equation}\label{eq-size-sifted}
  S(\mathcal{F},z) =|\{m\leq X\mid\, f(m)\text{ has no prime factor $<z$}\}|,
\end{equation}
and in particular, if $z\approx X^{1/2r}$, we get
$$
S(\mathcal{F},z) =|\{\text{primes $\gg X^{r/2}$ of the form $f(m)$
  with $m\leq X$}\}|,
$$
a function of much arithmetic interest.
\par
One immediately notices in this example that, to be effective, sieve
methods have to handle estimates uniform in terms of the support of
the sequence, and in terms of the parameter $z$ determining the set of
the primes in the sieve, as the latter will be a function of the
former.
\par
The congruence sums modulo $d\geq 1$ squarefree are easy to understand
here: we have
$$
S_d(\mathcal{F}_{f,X})=
\sum_{\stacksum{m\leq X}{d\mid f(m)}}{1},
$$
and by splitting the sum over $m$ into residue classes modulo $d$ and
denoting
$$
\rho_f(d)=|\{\alpha\in\Zz/d\Zz\,\mid\, f(\alpha)\equiv 0\mods{d}\}|
$$
the number of roots of $f$ modulo $d$, we find
\begin{equation}\label{eq-splitting}
S_d(\mathcal{F}_{f,X})=
\sum_{\stacksum{\alpha\in \Zz/d\Zz}{f(\alpha)=0\mods{d}}}{
\sum_{\stacksum{m\leq X}{m\equiv \alpha\mods{d}}}{1}
}=\frac{\rho_f(d)}{d}X+O(1)
\end{equation}
for $X\geq 2$, where the implied constant depends on $f$. The Chinese
Remainder Theorem shows that $d\mapsto \rho_f(d)$ is multiplicative,
and then from the Chebotarev density theorem (in the general case,
though a trivial argument suffices if $f(T)=(T-a_1)\cdots (T-a_r)$
splits completely over $\Zz$, which corresponds to the
Hardy-Littlewood prime tuple conjecture\footnote{\ Which played an
  important role in the recent results concerning small gaps between
  primes of Goldston, Pintz and Y\i ld\i r\i m.}), we know that
$$
\sum_{p\leq x}{\frac{\rho_f(p)}{p}}=\kappa \log\log x+O(1),
$$
where $\kappa=\kappa(f)$ is the number of irreducible factors of $f$
in $\Qq[T]$, so that the sieve problem here is of dimension
$\kappa$. (E.g., for $f(T)=T^2+1$, we have $\kappa=1$; there is, on
average over primes $p$, one square root of $-1$ in $\Zz/p\Zz$.)
\end{exem}

A number of very sophisticated combinatorial and number-theoretic
analysis (starting with Brun) have led to the following basic sieve
statement (see, e.g.,~\cite[Th. 11.13]{FI}, where more details are
given):

\begin{theo}\label{th-fund}
  With notation as above, for a sieve problem of dimension $\kappa>0$,
  there exists a real number $\beta(\kappa)>0$ such that
\begin{multline*}
(f(s)+O(\log D)^{-1/6}))S(\mathcal{F})\prod_{p\mid P(z)}{(1-g(p))}+R(D) \leq
S(\mathcal{F},P) \\
\leq (F(s)+O((\log D)^{-1/6}) S(\mathcal{F})\prod_{p\mid
  P(z)}{(1-g(p))}+R(D)
\end{multline*}
for $z=D^{1/s}$ with $s>\beta(\kappa)$, where $F(s)>0$ and $f(s)>0$
are certain functions of $s\geq 0$, depending on $\kappa$, defined as
solutions of explicit differential-difference equations, such that
$$
\lim_{s\ra +\infty}{f(s)}=\lim_{s\ra +\infty}{F(s)}=1,
$$
and where
$$
R(D)=\sum_{d<D}{|r_d(\mathcal{F})|}.
$$
\par
In both upper and lower bounds, the implied constant depends only on
$\kappa$ and on the constants in the
asymptotic~\emph{(\ref{eq-dimension2})}. 
\end{theo}

The leading term in the upper and lower bounds has a clear intuitive
meaning: one congruence condition leads to a proportion $1-g(p)$ of
the sum $S(\mathcal{F})$ that ``passes'' the test of sieving by $p$,
and multiple congruence conditions behave -- up to a point -- as if
they were independent. In particular, as soon as the remainders
in~(\ref{eq-sieve-axiom}) are small for fixed $d$, we have an
asymptotic formula for sieving with a \emph{fixed} set of primes, as
the support of the sequence $\mathcal{F}$ grows.
\par
Since $g(p)$ is about $\kappa p^{-1}$ on average, the Mertens formula
leads to
$$
\prod_{p\mid P}{(1-g(p))}\asymp \frac{1}{\log X}
$$
for $z=X^{1/s}$ for any fixed $s>0$. The following definition is
therefore a natural expression of the fact that one needs $R$ to be of
smaller order of magnitude to obtain actual consequences from
Theorem~\ref{th-fund}.

\begin{defi}[Level of distribution]\label{def-level}
  Let $(\mathcal{F}_n)$ be sequences as above and $D_n>0$. The
  $\mathcal{F}_n$ have \emph{level of distribution} $\geq D_n$ if
  and only if
\begin{equation}\label{eq-level}
R_n=\sum_{d<D_n}{|r_d(\mathcal{F}_n)|}\ll S(\mathcal{F}_n)(\log D_n)^{-B}
\end{equation}
for any $B>0$ and $n\geq 2$, the implied constant depending on $B$.
\end{defi}

\begin{exem}\label{ex-2}
  In the context of Example~\ref{ex-1} for $f\in\Zz[T]$ of degree $r$,
  with $\kappa$ irreducible factors, we have
$$
r_d(\mathcal{F}_{f,X})\ll d^{\eps}
$$
for all $d\geq 1$ squarefree and $\eps>0$, the implied constant
depending on $\eps$. Since $S(\mathcal{F}_{f,X})\asymp X$, we see that
the level of distribution is $\geq D$ for any $D=X^{1-\delta}$ for
$\delta>0$. Applying Theorem~\ref{th-fund} with $z=D^{1/s}$, $s$ large
enough, we deduce that there exists $\satur(f)\geq 1$ such that there
are infinitely many positive integers $m$ such that $f(m)$ has at most
$\satur(f)$ prime factors, counted with multiplicity (and in fact,
there are $\gg X/(\log X)^{\kappa}$ such integers $m\leq X$).
\end{exem}

\subsection{The sieve as a local-global study}
\label{ssec-local-global}

The point of view just described concerning sieves is very
efficient. However, we will use an essentially equivalent formal
description which is more immediately natural in the applications we
want to consider later.
\par
In this second viewpoint, we start with a set $Y$ of objects of
``global'' (and often arithmetic) nature. To study them, we assume
given maps
$$
Y\ra Y_p
$$
for $p$ prime which are analogues of (and often defined by)
``reduction modulo $p$''. Often, $Y$ parametrizes certain integers,
and the $Y_p$ parametrize their residue classes modulo $p$. To
emphasize this intuition, we write $y\mods{p}$ for the image in $Y_p$
of $y\in Y$. We interpret these maps as giving local information on
objects in $Y$, and we assume that $Y_p$ is a finite set. It is often
the case that $Y\ra Y_p$ is surjective, but sometimes it is convenient
to allow for non-surjective reduction maps.
\par
Using such data, we can form sifted sets associated with a set
$\mathcal{P}$ of primes, some $z\geq 2$, and some subsets
$\Omega_p\subset Y_p$, namely
\begin{equation}\label{eq-sifted}
\sifted_z(Y;\Omega)=
\{y\in Y\,\mid\, y\mods{p}\notin \Omega_p\text{ for all } 
p\in\mathcal{P},\ p<z\}\subset Y.
\end{equation}
\par
To ``count'' the elements in this sifted set, we consider quite
generally that we have available a finite measure $\mu$ on $Y$, and
the problem we turn to is to estimate (asymptotically, or from above
or below), the measure $\mu(\sifted_z(Y;\Omega))$ of the sifted
set.\footnote{\ It is hoped that $\mu$ will not be mistaken with the
  M\"obius function.} We assume for simplicity that for given $y\in
Y$, the set of $p\in\mathcal{P}$ with $y\mods{p}\in\Omega_p$ is finite
(this holds in most applications).
\par
This question can be interpreted in the previous framework as follows:
for any $y\in Y$, we define
$$
n(y)=\prod_{\stacksum{p\in\mathcal{P}}{y\mods{p}\in\Omega_p}}{p}
$$
for $y\in Y$, with the convention that $n(y)=0$ if the product is
infinite. This is a non-negative integer depending on $y$ such that the
``adjunction'' property
$$
(p\mid n(y))\Longleftrightarrow (y\mods{p}\in \Omega_p)
$$
holds for all $p\in\mathcal{P}$ if $n(y)\geq 1$. 
\par
Although the case $n(y)=0$ is usually exceptional,\footnote{\ In
  counting questions below, it will have a negligible contribution.}
it may occur. In order to take it into account, we define
$$
Y^0=\{y\in Y\,\mid\, n(y)=0\},\quad\quad Y^+=Y-Y^0.
$$
\par
We define the sequence $\mathcal{F}=(a_n)_{n\geq 1}$ by\footnote{\ We
  assume of course that all sets $\{y\,\mid\, n(y)=\alpha\}$ are
  measurable.}
\begin{equation}\label{eq-newa}
a_n=\mu(\{y\in Y\,\mid\, n(y)=n\}).
\end{equation}
\par
It follows that 
$$
S(\mathcal{F})=\sum_{n}{a_n}=\mu(Y^+),
$$
and
$$
S(\mathcal{F},z)= \sum_{(n,P(z))=1}{a_n}= \mu\Bigl( \{y\in Y^+\,\mid\,
(n(y),P(z))=1\} \Bigr)= \mu(\sifted_z(Y^+;\Omega)).
$$
\par
\begin{exem}\label{ex-1bis}
  Example~\ref{ex-1} may also be interpreted in this manner: we take
  $Y$ to be the set of positive integers, $\mu$ to be the counting
  measure restricted to integers $1\leq m\leq X$ in $Y$, the reduction
  maps to be $Y\ra \Zz/p\Zz$, which are indeed surjective maps onto
  finite sets. With
$$
\Omega_p=\{\alpha\in\Zz/p\Zz\,\mid\, f(\alpha)=0\}
$$
the set of zeros of $f$ modulo $p$, it is clear that
$$
\mu(\sifted_z(Y;\Omega))=S(\mathcal{F},z)
$$
is the quantity given in~(\ref{eq-size-sifted}).\footnote{\ The set
  $Y^+$ is here the set of $m\leq X$ such that $f(m)\not=0$; in
  particular, $Y^0$ is a finite set.}
\end{exem}

Similarly, the congruence sums $S_d(\mathcal{F})$ are now given by
$$
S_d(\mathcal{F})=\mu(\{y\in Y^+\,\mid\, y\mods{p}\in\Omega_p\text{ for
  all } p\mid d\})
$$
for $d$ squarefree. This is also the measure of the set
$$
\Omega_d=\prod_{p\mid d}{\Omega_p}\subset \prod_{p\mid d}{Y_p}
$$
under the image measure of $\mu$ by the map of simultaneous reduction
modulo $p\mid d$ on $Y^+$.
\par
It is therefore very natural to take the following point of view
towards the dimension condition~(\ref{eq-sieve-axiom}): first, we
expect that for $p$ prime -- provided the measure $\mu$ counts ``a
large part'' of $Y$ --, we can write
$$
\mu(\{y\in Y\,\mid\, y\mods{p}=\alpha\})\approx \mu(Y)\nu_p(\alpha)
$$
for all $\alpha$, where $\nu_p$ is some natural probability measure on
the finite set $Y_p$; secondly, we expect that for $d$ squarefree, the
reductions modulo the primes $p$ dividing $d$ are (approximately)
independent, so that
\begin{equation}\label{eq-equi-heuristique}
\mu(\{y\in Y\,\mid\, y\mods{p}=\alpha_p\text{ for all }
p\mid d\})\approx \mu(Y)\prod_{p\mid d}{\nu_p(\alpha_p)}.
\end{equation}
\par
If we compare this with~(\ref{eq-sieve-axiom}), this corresponds to
taking
$$
g(p)=\nu_p(\Omega_p),\quad\quad g(d)=\prod_{p\mid d}{\nu_p(\Omega_p)},
$$
and fits well with the intuitive meaning of $g(d)$ as encoding the
density of the sequence restricted to $n$ divisible by $d$.  Most
crucially, we note that assuming that $g$ is multiplicative is
essentially \emph{equivalent} with the expected asymptotic
independence property~(\ref{eq-equi-heuristique}). 
\par
We now proceed to make approximations~(\ref{eq-equi-heuristique})
precise, and interpret the level of distribution condition.  This
makes most sense in an asymptotic setting, and we therefore assume
that we have we have a sequence $(\mu_n)$ of finite measures on $Y$
(corresponding to taking $X\ra +\infty$ in Example~\ref{ex-1}).
\par
For any fixed squarefree number $p$, we consider the image
of the associated probability measures
$$
\tmu_n=\mu_n/\mu_n(Y),
$$
on the finite set
$$
Y_d=\prod_{p\mid d}{Y_d},
$$
which are probability measures $\tmu_{n,d}$ on $Y_d$.

\begin{defi}[Basic requirements]\label{def-require}
  Sieve with level $D_n\geq 1$ is possible for the objects in $Y$
  using $Y\ra Y_p$, and for a sequence $(\mu_n)$ of measures on $Y$,
  provided we have:
\par
\emph{(1) [Existence of independent local distribution]} For every $d$
squarefree, the image measures $\tmu_{n,d}$ converge to some
measure $\nu_d$, and in fact
$$
\nu_d=\prod_{p\mid d}{\nu_p}.
$$
\par
In other words, there are probability measures $\nu_d$ on $Y_d$ such
that
$$
\mu_n(\{y\in Y\,\mid\, y\mods{p}=\alpha_p\})
=\mu_n(Y)(\nu_d(\alpha)+r_{d,n}(\alpha))
$$
for all $d$ and $\alpha=(\alpha_p)_{p\mid d}\in Y_d$, and
$$
\lim_{n\ra +\infty}{r_{d,n}(\alpha)}=0
$$
for all fixed $d$ and $\alpha\in Y_d$.
\par
\emph{(2) [Level of distribution condition]} Given subsets
$\Omega_p\subset Y_p$, and
$$
\Omega_d=\prod_{p\mid d}{\Omega_p},
$$
we have
\begin{gather}\label{eq-level-bis}
\sum_{d<D_n}{|\Omega_d|\max_{\alpha}|r_{d,n}(\alpha)|}
\ll (\log D_n)^{-A}\\
\label{eq-fin-unite}
  \mu_n(Y^0)\ll \mu_n(Y)(\log D_n)^{-A}
\end{gather}
for all $n\geq 1$, the implied constant depending on $A$ and the
$\Omega_p$. 
\end{defi}

These assumptions are a form of \emph{quantitative}, \emph{uniform}
equidistribution results for the reductions of objects of $Y$ modulo
primes, and of independence of these reductions modulo various
primes. Much of the difficulty of the sieve in orbits (as in
Section~\ref{ssec-apol} or more generally, as described in the next
section) lies in checking the validity of these conditions. We will
present some of the techniques and new results in the next
sections. 
\par
However, if we assume that the basic requirements are met, we may
combine these conditions with the fundamental statement of sieve, to
deduce the following general fact. Given sets $\Omega_p\subset Y_p$
such that
\begin{equation}\label{eq-dim-kappa}
\nu_p(\Omega_p)=\frac{\kappa}{p}+O(p^{-1-\delta})
\end{equation}
(a condition which is local and does not depend at all on $Y$) for
some $\delta>0$ and $p$ prime, we obtain by summing over $\Omega_p$
that the congruence sums satisfy
$$
\mu_n(y\in Y\,\mid\, y\mods{p}\in\Omega_p\text{ for } p\mid d)=
\mu_n(Y)\nu_d(\Omega_d)+r_{d,n}(\Omega)
$$
where 
$$
r_{d,n}(\Omega)\ll |\Omega_d|\max_{\alpha\in Y_d}{|r_{d,n}(\alpha)|}.
$$
\par
By independence, $g(d)=|\Omega_d|$ is multiplicative. Hence we derive
(weakening the conclusion of Theorem~\ref{th-fund}) that
\begin{equation}\label{eq-easy-sieve}
\mu_n(\{y\in Y\,\mid\, y\mods{p}\notin\Omega_p\text{ for all }
p<D_n^{1/s}\})\asymp
\frac{\mu_n(Y)}{(\log D_n)^{\kappa}}
\end{equation}
for $n\geq 1$ and for all fixed $s$ large enough (in terms of
$\kappa$).

\begin{exem}
  This description applies easily to Examples~\ref{ex-1}
  and~\ref{ex-1bis} with the sequence $(\mu_X)$ of counting measures
  on $1\leq m\leq X$. The measure $\nu_p$ is simply the uniform
  probability measure on $\Zz/p\Zz$, and the independence (which
  passed almost unnoticed) is valid, as an expression of the Chinese
  Remainder Theorem: knowing the reduction modulo a prime $p_1$ of a
  general integer $n$ gives no information whatsoever on its reduction
  modulo another prime $p_2$. Quantitatively, we have 
$$
|\{m\leq X\,\mid\, m\equiv \alpha\mods{d}\}|= \frac{X}{d}+O(1),
$$
and we recover the level of distribution $D=X^{1-\delta}$ with
$\delta>0$ of Example~\ref{ex-2}. (This is about the best one can hope
for.) Note that condition~(\ref{eq-fin-unite}) is here trivial since
$Y^0$ is the finite set of positive integral roots of $f$.
\end{exem}


\begin{rema}\label{rm-finitude}
  The upper-bound~(\ref{eq-fin-unite}) is obtained, in all the cases
  discussed in this report, as an immediate consequence
  of~(\ref{eq-level-bis}). Indeed, in all cases we consider, it will
  be true that $n(y)=0$ means that $y\mods{p}\in\Omega_p$ for
  \emph{all} $p$. One can then write
$$
\mu_n(Y^0)\leq \mu_n(\{y\in Y\,\mid\, y\mods{p}\in\Omega\})
$$
for any fixed $p$, and under the assumptions above, we find for $p\leq
D_n$ that 
\begin{align*}
\mu_n(Y^0)&\leq \mu_n(Y)\Bigl\{\nu_p(\Omega_p)+O(|\Omega_p|\max_{\alpha\in
  \Omega_p}{|r_{p,n}(\alpha)|})\Bigr\}
\\
&\ll \mu_n(Y)\Bigl\{\frac{1}{p}+O((\log D_n)^{-A})\Bigr\},
\end{align*}
so we derive~(\ref{eq-fin-unite}) by taking $p$ of size comparable
with $(\log D_n)^A$.
\end{rema}

\section{Sieve in orbits}\label{sec-orbits}

We will now present the general version of the sieve problem developed
by Bourgain, Gamburd and Sarnak, which generalizes the example of
Apollonian circle packings (Section~\ref{ssec-apol}).  As with the
group $\apol$, the global objects of interest are directly related to
discrete groups with exponential growth. (Further discussion of the
example in Section~\ref{ssec-geom} is found in
Section~\ref{sec-geom}.)

\subsection{The general setting}

Consider a finitely-generated group $\Lambda\subset \GL_m(\Zz)$ for
some $m\geq 1$ (e.g., $\apol\subset \GL_4(\Zz)$, the ``Lubotzky''
group $L$ of~(\ref{eq-lub}) in $\SL_2(\Zz)$, or $\SL_m(\Zz)$ itself).
\par
Given a non-zero vector $x_0\in \Zz^m$, form the $\Lambda$-orbit
$$
\orbit(x_0)=\Lambda \cdot x_0\subset \Zz^m,
$$ 
and fix some polynomial function $f\in \Qq[X_1,\ldots,X_m]$ such that
$f$ is integral valued and non-constant on $\orbit(x_0)$ (for
instance, $f\in \Zz[X_1,\ldots, X_m]$). The ``philosophical'' question
to consider is then
\par
\begin{center}
``To what extent are the values $f(x)$, where $x\in
\orbit(x_0)$, typical integers?''
\end{center}
\par
More precisely, as far as sieve is concerned, the question is: how do
the multiplicative properties (number and distribution of prime
factors) of the $f(\gamma x_0)$ differ from those of general
integers?\footnote{\ For readers unfamiliar with what this means, the
  Appendix gives some of the most basic statements along these lines.}
From the point of view of Section~\ref{ssec-local-global}, we are
therefore looking at $Y=\orbit(x_0)$, and its reductions modulo primes
$$
Y\ra Y_p=\orbit(x_0\mods{p}),
$$
the latter being the orbit of $x_0\mods{p}$ under the action of
$\Lambda_p$, the image of $\Lambda$ under reduction modulo primes.

\begin{rema}
  One may also consider other actions of $\Lambda$; for instance
  taking $Y=\Lambda$ itself is natural enough, with $Y_p$ the
  reduction of $\Lambda$ modulo $p$. However, if one considers the
  image of $\Lambda$ in $\GL_{m^2}(\Zz)$ corresponding to the
  multiplication action of $\Lambda$ on $\GL_m(\Zz)$ on the left, the
  group $\Lambda$ is isomorphic to the orbit of the vector
  ``identity'' $x_0=1\in M_m(\Zz)\simeq \Zz^{m^2}$.
\end{rema}

Later we will explain how to count elements of $Y$ in order to apply
the basic sieve statements. However, there is a first qualitative way
of phrasing the guess that there should be many elements $x$ of
$\orbit(x_0)$ with $f(x)$ having few prime factors, which was pointed
out in~\cite{bgs}. Namely, let 
$$
\orbit_f(x_0;\satur)= \{ x\in\orbit(x_0)\,\mid\, \Omega(f(x))\leq
r \}
$$
for $r\geq 1$ and define the ``saturation number'' for the orbit:
$$
\satur(f,\Lambda)=
\min\{r\geq 1\,\mid\, \orbit_f(x_0;r)\text{ and }
\orbit(x_0)\text{ have the same Zariski-closure}\},
$$
or in other words, the smallest $r\geq 1$ such that the elements of
$\orbit_f(x_0;r)$ satisfy no further polynomial relation than those of
the full orbit $\orbit(x_0)$. One asks then:

\begin{quest}\label{quest-satur}
Is this ``saturation number'' finite, and if yes, what is it? 
\end{quest}

\begin{exem}\label{rm-pythagor}
  When $m=1$, a set in $\Zz$ is either Zariski-closed, if finite, or
  has Zariski-closure the whole affine line. Thus the existence of
  saturation number for an infinite subset $\orbit$ of $\Zz$ amounts
  to no more (but no less) than the statement that $\orbit$ is
  infinite. 
\par
However, with $m\geq 2$, the saturation condition may become quite
interesting. The following example is discussed in~\cite[\S
6,Ex. C]{bgs}: consider the orbit $\orbit$ of $x_0=(3,4,5)$ under the
action of the orthogonal group $\Lambda=\SO(2,1)(\Zz)$. This orbit is
the set of integral Pythagorean triples, and its Zariski closure is
the cone $\{x^2+y^2-z^2=0\}$. Considering the function $f(x,y,z)=xy/2$
(the area of the right-triangle associated with $(x,y,z)$), and using
the recent \emph{quantitative} results of Green and Tao~\cite{gt} on
the number of arithmetic progressions of length $4$ in the primes,
Bourgain, Gamburd and Sarnak show that the saturation number is $6$ in
that case. However, it is highly likely (this follows from some of the
Hardy-Littlewood conjectures) that there are infinitely many integral
right-triangles with area having $\leq 5$ prime factors! These
triangles, however, have sides related by a non-trivial polynomial
relation.
\end{exem}

It seems interesting to to sharpen a bit the strength of the
saturation condition by replacing the condition ``Zariski-dense'' with
the stronger condition that $\orbit_f(x_0;r)$ be \emph{not thin}. We
recall the definition (see~\cite[\S 3.1]{serre-galois}):

\begin{defi}[Thin set]\label{def-thin}
  Let $V/k$ be an irreducible algebraic variety defined over a field
  $k$ of characteristic zero. A subset $A\subset V(k)$ is \emph{thin}
  if there exists a $k$-morphism $W\fleche{f} V$ of algebraic
  varieties such that $f$ has no $k$-rational section, $\dim (W)\leq
  \dim(V)$ and
$$
A\subset f(W(k)).
$$
\end{defi}

\begin{exem}
For $m=1$, there are many infinite thin sets in $\Zz$. However, one
shows that such a set (say $\mathcal{T}$) satisfies
$$
|\{n\in \mathcal{T}\,\mid\, |n|\leq X\}|\ll X^{1/2}(\log X)
$$
for $X\geq 2$ (a result of S.D. Cohen, see,
e.g.,~\cite[Th. 3.4.4]{serre-galois}), so that even Chebychev's
elementary bounds prove that the set of prime numbers is not thin. The
example of the set of squares shows that the exponent $1/2$ is best
possible.
\end{exem}

We now see that Theorem~\ref{th-bgs-ex} states, for certain groups
$\Lambda$ and their orbits, that the saturation number is finite, even
that $\orbit_f(x_0;r)$ is non-thin. We will sketch below the proof of
this result, in a slightly more general case.  Interestingly,
although~\cite{bgs} only considers the original saturation condition,
the values of $r$ for which they prove that $\orbit_f(x_0;r)$ is
Zariski-dense are always such that it is not thin.


\subsection{How to count?}
\label{sec-counting}

In the setting of the previous section, we have many examples of a set
$Y$ with reduction maps $Y\ra Y_p$, and we wish to implement the sieve
techniques, for instance to prove that some saturation number is
finite. For this purpose, as described in Section~\ref{sec-sieve}, it
is first necessary to specify how one wishes to count the elements of
$Y$, or in other words, one must specify which finite measures
$(\mu_n)$, will be used to check the Basic Requirements of sieve
(Definition~\ref{def-require}).
\par
A beautiful, characteristic, feature of the sieve in expansion is that
there are often two or three natural ways of counting the elements of
$Y$. We illustrate this in the case that $Y=\Lambda$ is a (finitely
generated) subgroup of $\GL_m(\Zz)$. One may then use:
\par
-- [Archimedean balls] One may fix a norm $\|\cdot\|$ on $\GL_m(\Rr)$,
and let $\mu_X$ be the uniform counting measure on the finite set
$$
B_{\Lambda}(X)=\{g\in \Lambda\subset \GL_m(\Rr)\,\mid\, \|g\|\leq X\}
$$
for some $X\geq 1$. 
\par
-- [Combinatorial balls] One may fix instead a finite generating set
$S$ of $\Lambda$, assumed to be symmetric (i.e., $s\in S$ implies
$s^{-1}\in S$), and use it to define first a combinatorial word-length
metric, i.e.,
$$
\|g\|_S=\min\{k\geq 0\,\mid\, g=s_1\cdots s_k\text{ for some }
s_1,\ldots, s_k\in S\}.
$$
\par
Then one can then consider the uniform counting probability measure
$\mu_k$ on the finite combinatorial ball
$$
B_S(k)=\{g\in \Lambda\,\mid\, \|g\|_S\leq k\}
$$
for $k\geq 1$ integer. This set depends on $S$, but many robust
properties should be (and are) independent of the choice of generating
sets.
\par
-- [Random walks] Instead of the uniform probability on combinatorial
balls, it may be quite convenient to use a suitable weight on the
elements of $B_S(k)$ that takes into account the multiplicity of their
representation as words of length $k$. More precisely, we assume that
$1\in S$ (adding it up if necessary) and we consider the probability
measure $\mu_k$ on $\Lambda$ such that
$$
\mu_k(g)=\frac{1}{|S|^k} \multsum_{\stacksum{(s_1,\ldots, s_k)\in
    S^k}{s_1\cdots s_k=g}}{ 1 }
$$
(adding $1$ to a generated set $S$, if needed, ensures that this
measure is supported on the $S$-combinatorial ball or radius $k$, and
not the combinatorial sphere.)
\par
Compared with the previous combinatorial balls, the point of this
weight is that it simplifies enormously any sum over $B_S(k)$: we have
$$
\sum_{g\in \Lambda}{\varphi(g)\mu_k(g)}=
\frac{1}{|S|^k}
\multsum_{s_1,\ldots, s_k\in S}{
\varphi(s_1\cdots s_k)
},
$$
for any function $\varphi$ on $\Lambda$, where the summation variables
on the right are \emph{free}.

\begin{rema}
  This third weight has a natural probabilistic interpretation:
  $\mu_k$ is the law of the $k$-th step of the left-invariant random
  walk $(X_k)_{k\geq 0}$ on $\Lambda$, defined by
$$
X_0=1\in \Lambda\,\quad\, X_{k+1}=X_k\xi_{k+1},
$$
where $(\xi_k)_{k\geq 1}$ is a sequence of $S$-valued independent
random variables (on some probability space $(\Omega,\Sigma,\proba)$)
such that
$$
\proba(\xi_k=s)=\frac{1}{|S|}\quad\text{ for all $k\geq 1$ and $s\in S$.}
$$
\end{rema}

Once the counting method $(\mu_X)$ is chosen, the more precise form of
Question~\ref{quest-satur} becomes to bound from below the function
$$
\pi_f(X;r)=\mu_X(\gamma\in \Lambda\,\mid\, \Omega(f(\gamma\cdot
x_0))\leq r)
$$
as $X\ra +\infty$. Precisely, the idea is to prove -- for suitable $r$
-- a lower bound which is sufficient to ensure that $\orbit_f(x_0;r)$
is Zariski-dense, by comparison with upper bounds (known or to be
established) for the counting functions
$$
\mu_X(\gamma\in \Lambda\,\mid\, f(\gamma\cdot x_0)\in W)
$$
for a proper Zariski-closed subvariety $W\subset
V=\overline{\orbit(x_0)}$ (or for a thin subset $W\subset V(\Qq)\cap
\Zz^m$).

\section{The basic requirements}

Consider $\Lambda$ and an orbit $\orbit(x_0)$ as in the previous
section, and assume a counting method (i.e., a sequence $(\mu_X)$ of
measures on $Y$) has been selected. We now proceed to check if the
basic requirements of sieve are satisfied.
\par
However, we first make the following assumption (which will be refined
later):

\begin{hypot}
  The Zariski-closure $G/\Qq$ of the group $\Lambda\subset \GL_m(\Zz)$
  is a \emph{semisimple group}, e.g., $\SL_m$, or $\Sp_{2g}$ or an
  orthogonal group, or a product of such groups.
\end{hypot}

For instance, this means that we exclude from the considerations below
the subgroup generated by the single element $2\in \GL_1(\Zz[1/2])$
and its orbit $\{2^n\}\subset \Zz[1/2]$ (we extend here the base ring
slightly); this is understandable because the question of finding
integers $n$ for which, say, $f(2^n)=2^n-1$, is prime remains
stubbornly resistant. Indeed, the basic results used to understand the
local information available from $\Lambda\ra \Lambda_p$ simply fail in
that case (see Remark~\ref{rm-mersenne} below). In~\cite[\S 2]{bgs},
Bourgain, Gamburd and Sarnak give some more examples that show why
general reductive groups lead to very different -- and badly
understood -- phenomena.
\par
The same reason make solvable groups (e.g., upper-triangular matrices)
delicate to handle; as for nilpotent groups, they are definitely
accessible to sieve methods. However, their behavior (in terms of
growth) is milder and one would not need the considerations of
expansion in groups that are needed for semisimple groups.

\subsection{Local limit measures}

According to Section~\ref{ssec-local-global}, we start the
investigation of sieve by checking whether, for a fixed squarefree
integer $d$, the measures $\tmu_{X,d}$ on
$$
\prod_{p\mid d}{\Lambda_p}
$$
have a limit, and whether this limit is a product measure.  This last
condition is crucial, and it sometimes require more preliminary
footwork. The basic difficulty is illustrated in the following
situation: suppose that, for some set $Z$, there are non-constant maps
$$
N\,:\, Y\ra Z,\quad\quad N_p\,:\, Y_p\ra Z,
$$
and that
$$
N(y)=N_p(y\mods{p})
$$
for all primes $p$. Then $y\mods{p}$ \emph{does} carry some
information concerning $y$, namely the value of $N(y)$. Consequently,
it is not possible for $y\mods{p_1p_2}$ to become equidistributed
towards a product measure $\nu_{p_1}\times \nu_{p_2}$.

\begin{exem}
In the sieve in orbits, this happens frequently when orthogonal groups
are concerned. For instance, the Apollonian group $\apol$ is not
contained in $\SO(Q,\Zz)$ and we have
$$
\det(\gamma)=\det(\gamma\mods{p})
$$
for all $p$. (In fact, there is a further obstruction even for
$\SO(Q,\Zz)$.)
\end{exem} 

If such obstacles to independence appear for a given problem, this is
a sign that it should reworded or rearranged: instead of $Y$, one
should attempt to sieve, for instance, the fibers of $N$. (This is
justified by the fact that, in some cases, different fibers may indeed
have very different behavior, and can not be treated uniformly).
\par
In our case of sieve in orbits for a group $\Lambda$ which is
Zariski-dense in the semisimple group $G/\Qq$, the most efficient
method is to reduce first to the connected component $G^0$ of $G$ (by
replacing $\Lambda$ with $\Lambda\cap G^0(\Qq)$ and then to the
\emph{simply-connected} covering group $G^{sc}$ of $G^0$ using the
projection map
$$
\pi\,:\,G^{sc}\ra G^0\ ;
$$
one works with the inverse image $\Lambda^{sc}$ of $\Lambda$ in
$G^{sc}(\Qq)$, and consider the function $\tilde{f}=f\circ\pi$ instead
of $f$. Note that, in general, either of these operations might
require to work over a base number field distinct from $\Qq$, but
there is no particular difficulty in doing so. For a detailed analysis
in the case of the Apollonian group $\apol$, where $G=O(Q)$ is not
connected and the connected component $\SO(Q)$ is not
simply-connected, see~\cite[\S 6]{bgs} or~\cite{fuchs}.

\begin{exem}
  In the situation of Theorem~\ref{th-bgs-ex}, we have $G=\SL_m$,
  which is connected and simply-connected, so these preliminaries are
  not needed. The same thing happens if $G$ is a symplectic group
  $G=\Sp_{2g}$, or if $G$ is a product of groups of these two types.
\end{exem}

The following result is the crucial ingredient that shows that a
Zariski-dense subgroup in a simply-connected group satisfies a strong
form of independence of reduction modulo primes:

\begin{theo}[Strong approximation and independence]\label{th-approx}
  Let $G$ be a connected, simply-connected, absolutely almost simple
  $\Qq$-group embedded in $\GL_m/\Qq$, and let $\Lambda\subset
  G(\Qq)\cap \GL_m(\Zz)$ be a Zariski-dense subgroup.\footnote{\ For
    instance $G=\SL_m$ or $\Sp_{2g}$.} There exists a finite set of
  primes $\sop=\sop(\Lambda)$ such that $G$ has a model, still denoted
  $G$, over $\Zz[1/\sop]$, and:
\par
\emph{(1)} For all primes $p$ not in $\sop$, the map
$$
\Lambda\ra G(\Fp_p)
$$
is surjective, i.e., the image $\Lambda_p$ of reduction modulo $p$ is
``as large as possible'', so that $\Lambda_p=G(\Fp_p)$.
\par
\emph{(2)} For all squarefree integers $d$ coprime with $\sop$, the
reduction map
$$
\Lambda \ra \prod_{p\mid d}{G(\Fp_p)}=G(\Zz/d\Zz)
$$
is surjective, i.e., we have $\Lambda_d=G(\Zz/d\Zz)$ and $\Lambda\ra
\Lambda_d$ is surjective.
\par
\emph{(3)} Assume $\mu_k$ is the weighted counting method of
Section~\ref{sec-counting} associated with a finite symmetric
generating set $S$, with $1\in S$. Then, for any integer $d$ coprime
with $\sop$, the probability measures $\tmu_{k,d}$ on
$$
\Lambda_d=\prod_{p\mid d}{\Lambda_p}=G(\Zz/d\Zz)
$$
converge, as $k\ra +\infty$, to the uniform probability measure
$\nu_d=\prod{\nu_p}$, i.e., the measure so that
$$
\nu_d(\gamma)=\frac{1}{|G(\Zz/d\Zz)|},\quad\quad
\text{ for all } \gamma\in G(\Zz/d\Zz).
$$
\end{theo}

Parts (1) and (2) have been proved, in varying degree of generality
(and with very different methods), by a number of people, in
particular Hrushovski and Pillai~\cite{hru-p}, Nori~\cite{nori},
Matthews-Vaserstein-Weisfeiler~\cite{mvw}; the most general statement
is due to Weisfeiler~\cite{weisfeiler}. 
\par
Part (3), on the other hand, is usually not proved a priori, It holds,
in fact, also in many cases when the two other counting methods are
used (i.e., archimedean balls and unweighted combinatorial balls) but
it is then seen as a consequence the stronger quantitative forms of
equidistribution (which are parts of the Basic Requirements
anyway). We will say more about this in the next sections, but we
should observe, However, that these limiting measures are certainly
the most natural ones that one might expect (being the Haar
probability measures on the finite groups $G(\Zz/d\Zz)$).
\par
We explain the quite straightforward proof of (3) for the weighted
counting (it is also an immediate consequence of the probabilistic
interpretation and standard Markov-chain theory).
\par
Let $\varphi\,:\, \Lambda_d\ra \Cc$ be any function. The integral of
$\varphi$ according to $\tmu_{k,d}$ is 
$$
\sum_{y\in\Lambda_d}{\varphi(y)
\frac{1}{|S|^k}
\multsum_{\stacksum{s_1,\ldots, s_k\in S}{s_1\cdots s_k=y}}{
1}}=
\frac{1}{|S|^k}\multsum_{s_1,\ldots, s_k\in S}{\varphi(s_1\cdots s_k)}
=(M^k\varphi)(1)
$$
where $M$ is the Markov averaging operator on functions on $\Lambda_d$
associated to $S$, i.e., for any $f\,:\, \Lambda_d\ra \Cc$ and
$x\in\Lambda_d$, we have 
$$
(Mf)(x)=\frac{1}{|S|}\sum_{s\in S}{f(xs)}.
$$
\par
The constant function $1$ is an eigenfunction of $M$ with eigenvalue
$1$. Because $S$ generates $\Lambda$ (hence $\Lambda_d$), it is an
eigenvalue with multiplicity $1$. Thus if we write
$$
\varphi=\frac{1}{|\Lambda_d|}\sum_{y\in
  \Lambda_d}{\varphi(y)}+\varphi_0,
$$
we have
$$
M^k\varphi=\frac{1}{|\Lambda_d|}\sum_{y\in
  \Lambda_d}{\varphi(y)}+M^k\varphi_0,
$$
and therefore
\begin{equation}\label{eq-equi-bound}
  \Bigl| (M^k\varphi)(1)-\frac{1}{|\Lambda_d|}\sum_{y\in
    \Lambda_d}{\varphi(y)}\Bigr| \leq \max_{y\in\Lambda_d}
  |(M^k\varphi_0)(y)|
  \leq |\Lambda_d|\rho_0(M)^k\|\varphi\|_2
\end{equation}
where $\rho_0$ is the spectral radius of the operator $M$ restricted
to the space of functions on $\Lambda_d$ with mean $0$ (according to
the uniform probability measure), endowed with the corresponding
$L^2$-norm $\|\cdot\|_2$. Having ensured that $1\in S$, as we required
for the weighted counting method, also implies that $-1$ is
not\footnote{\ Let $S'=S-\{1\}$, $|S|=s\geq 1$; the operator $M$ can
  be written $(1-1/s)M'+1/s$, $M'$ being the averaging operator for
  the generators $S'$; since the spectrum of $M'$ lies in $[-1,1]$,
  that of $M$ must be in $[-1+2/s,1]$.} an eigenvalue of $M$. Since
$M$ is symmetric, hence has real spectrum, and this spectrum is
clearly in $[-1,1]$, it follows that $\rho_0(M)<1$. Thus
$$
\int_{\Lambda_d}{\varphi(y)d\tmu_{k,d}}
\ra \frac{1}{|\Lambda_d|}\sum_{y\in
  \Lambda_d}{\varphi(y)}
$$
as $k\ra \infty$, which is the desired local equidistribution with
respect to $\nu_d$.

\begin{rema}\label{rm-mersenne}
  In fact, independence, in the sense of Theorem~\ref{th-approx}, only
  holds for simply-connected groups. Thus, its conclusion may be taken
  as a ``practical'' alternate characterization, with a very obvious
  interpretation, as far as sieve is concerned at least.
\par 
Now, consider again the cyclic group generated by $2\in
\GL_1(\Zz[1/2])$. Its image in $\GL_1(\Fp_p)$ (for $p$ odd) is cyclic,
and of order the order of $2$ modulo $p$. This remains a very
mysterious quantity; in particular, it is certainly not the case that
$2$ generates $\Fp_p^{\times}$ for most $p$ (a well-known conjecture
of Artin states that this should happen for a positive proportion of
the primes, but even this remains unknown). Thus the basic principles
of sieve break down at this early stage for this group. (The best that
has been done, for the moment, is to use the fact that, for almost all
primes $\ell$, there is some prime $p$ for which the order of $2$
modulo $p$ is $\ell$, in order to show that $2^n-1$ has, typically,
roughly as many \emph{small} prime factors as one would expect in view
of its size; see~\cite[Exercise 4.2]{cup}).
\end{rema}

\begin{exem}
For the Apollonian group, Fuchs~\cite{fuchs} has determined explicitly
the image under reduction modulo $d$ of (the inverse image in the
simply connected covering of) $\apol$, for all $d$.
\par
For concrete groups, it might also be possible to check Strong
Approximation directly. For the group $\lubotz\subset \SL_2(\Zz)$
(see~(\ref{eq-lub})) which has infinite index, for instance, it is
clear that $\lubotz$ surjects to $\SL_2(\Fp_p)$ for all primes
$p\not=3$, and that it is trivial modulo $3$. Then one obtains the
surjectivity of
$$
\lubotz\ra \prod_{p\mid d}{\SL_2(\Fp_p)}
$$
for all $d$ with $3\nmid d$ using the Goursat Lemma of group theory
(which is crucial in the proof of (2) in any case).
\end{exem}

\subsection{Quantitative equidistribution: combinatorial
  aspects}\label{sec-expanders}

We consider a finitely generated group $\Lambda\subset \GL_m(\Zz)$
such that its Zariski-closure $G$ is simple, connected and simply
connected, together with a symmetric finite generating set $S$. 
\par
If we select the weighted combinatorial count (assuming that $1\in
S$), we see from Theorem~\ref{th-approx} that in order to satisfy the
basic requirements of Definition~\ref{def-require} for sieving up to
some level $D_k$, we ``only'' must check the level of distribution
condition~(\ref{eq-level-bis}), provided we restrict our attention to
sieving using primes $p$ outside the possible finite exceptional set
$\sop$, and squarefree integers $d$ coprime with $\sop$.
\par
In turn, the estimate~(\ref{eq-equi-bound}), with $\varphi$ the
characteristic function of a point $\alpha$, shows that we have a
quantitative equidistribution result for fixed $d$, $(d,\sop)=1$, of
the type
$$
|r_{d,k}(\alpha)|\leq  \sqrt{|\Lambda_d|}\rho_d^k,
$$
valid for all $\alpha\in \Lambda_d$, where $\rho_d$ is the spectral
radius of the averaging operator on functions of mean zero on
$\Lambda_d$, which satisfies $\rho_d<1$.
\par
It is therefore obvious that \emph{obtaining a level of distribution
  for $\Lambda$ amounts to having an upper bound for $\rho_d$ which is
  uniform in $d$.} 
\par
The best that can be hoped for is that $\rho_d\leq \rho <1$ for all
$d\geq 1$ and some fixed $\rho$: the exponential rate of equidistribution
is then \emph{uniform} over $d$. This is equivalent to the well-known,
incredibly useful, condition that the family of Cayley graphs of
$\Lambda_d$ with respect to $S$ be an \emph{expander} (see~\cite{hlw}
for background on expanders). Precisely:

\begin{defi}[Expander family, random-walk definition]\label{def-expand}
  Let $(\Gamma_i)_{i\in I}$ be a family of connected $k$-regular
  graphs for a fixed $k\geq 1$, possibly with multiple edges or
  loops. The family $(\Gamma_i)$ is an expander if there exists
  $\rho<1$, independent of $i$, such that
$$
\rho_i\leq \rho<1
$$
for all $i$, where $\rho_i$ is the spectral radius of the Markov
averaging operator on functions of mean $0$ on $\Gamma_i$, with
respect to the inner product
$$
\langle f_1,f_2\rangle=\frac{1}{|\Gamma_i|}\sum_{x\in\Gamma_i}{
f_1(x)\overline{f_2(x)}
}.
$$
\end{defi}

Assuming that we have an expander, with expansion constant $\rho$, and
that
$$
|\Omega_p|\leq p^{\Delta},\quad |\Lambda_p|\leq p^{\Delta_1},
$$
we obtain the immediate estimate
$$
\sum_{d<D}{|\Omega_d|
\max_{\alpha\in\Lambda_d}|r_{d,k}(\alpha)|} \leq D^{1+\Delta+\Delta_1/2}\rho^k
$$
for $k\geq 1$. This provides levels of distribution (as in
Definition~\ref{def-require}) of the type
\begin{equation}\label{eq-delta}
  D_k=\beta^k,\quad \text{ for any }\quad 1<\beta<\rho^{-1/(1+\Delta+\Delta_1/2)}.
\end{equation}
\par
This applied to the weighted counting. One may however expect that,
under the condition that the Cayley graphs be expanders, a similar
level of distribution should hold for combinatorial balls
also. However, to the author's knowledge, this is currently only known
under the additional condition that $\Lambda$ be a (necessarily
non-abelian) free group. When this is the case, Bourgain, Gamburd and
Sarnak~\cite[\S 3.3, (3.29), (3.32)]{bgs} prove, using the well-known
spectral theory on free groups, that for the counting measure on the
combinatorial ball of radius $k$. there exists $\tau<1$, depending
only on the expansion constant $\rho$ of the family of graphs, such
that the measures $\tmu_{d,k}$ converge towards the uniform measure
$\nu_d$ on $\Lambda_d$ (for $d$ coprime to a suitable finite set of
primes), with error bounded by
$$
|r_{d,k}(\alpha)|\ll  |B_S(k)|^{\tau-1}
$$
for all $k\geq 1$.
\par
The restriction to free groups is awkward from some points of
view. However, for instance, the inverse image of $\apol$ in the
simply-connected covering of $\SO(Q)$ is free, so this can be used in
the case of the Apollonian circle packings. Moreover, as observed
in~\cite{bgs}, for applications such as upper-bounds on saturation
numbers, one can use the fact (the ``Tits alternative'') that, under
our assumptions on $G$, any Zariski-dense subgroup $\Lambda$ of $G$
contains a free subgroup $\Lambda_Z$ which is still Zariski-dense in
$G$. It is then possible to apply sieve to the orbit of $x_0$ under
$\Lambda_Z$ instead of $\Lambda$.
\par
It becomes, in any case, of pressing concern to know whether the
expansion property holds. Here we can reformulate the question in
terms of the Cayley graphs of $G(\Zz/d\Zz)$ for $d\geq 1$ squarefree,
with respect to suitable sets of generators (avoiding a few
exceptional primes).
\par
This question has in fact quite a pedigree. But until very recently,
the known examples where all lattices in semisimple groups.\footnote{\
  Except for some isolated examples of Shalom~\cite[Th. 5.2]{shalom}
  and, rather implicitly, the examples arising from Gamburd's
  work~\cite{gamburd} on the ``spectral'' side, detailed in the next
  section. (Both applied only to $d=p$ prime.)} Indeed, the expansion
property can also be phrased in the case of interest as stating that
the group $\Lambda$ has Property $(\tau)$ of Lubotzky for
representations factoring through the congruence quotients $\Lambda\ra
\Lambda_d=G(\Zz/d\Zz)$. Thus it is known, by work of
Clozel~\cite{clozel}, for $\Lambda$ any finite-index subgroup of
$G(\Zz)$. In fact, for many cases of great interest, like (finite
index subgroups in) $\SL_m(\Zz)$ for $m\geq 3$ or $\Sp_{2g}(\Zz)$ for
$g\geq 2$, the result follows from Kazhdan's Property (T).
\par
For subgroups $\Lambda$ (possibly) of infinite index in $G(\Zz)$,
there has been dramatic progress recently,\footnote{\ A story which is
  well worth its own account; the excellent survey~\cite{green} by
  B. Green, despite being very recent, does not cover many of the most
  remarkable new results.}  partly motivated by the sieve
applications. The following theorem has now been announced:

\begin{theo}[Expansion in finite linear groups]\label{th-expansion-groups}
  Let $G/\Qq$ be absolutely almost simple, connected and
  simply-connected, embedded in $\GL_m$ for some $m$, e.g., $G=\SL_m$,
  $m\geq 2$, or $\Sp_{2g}$, $g\geq 1$. Let $\Lambda\subset G(\Qq)\cap
  \GL_m(\Zz)$ be a finitely-generated subgroup which is Zariski-dense
  in $G$. Fix a symmetric finite system of generators $S$ of
  $\Lambda$. Then the family of Cayley graphs, with respect to $S$, of
  the groups $\Lambda_d$ obtained by reduction modulo $d$ of $\Lambda$
  is an expander family, where $d$ runs over squarefree integers.
\end{theo}

This applies, in particular, to the group $L$ of~(\ref{eq-lub}); this
was essentially a question of Lubotzky.
\par
Many people have contributed (and still contribute) to the proof of
this result (and variants, extensions, etc). We do not attempt a
complete history or any semblance of proof, but it seems useful to
sketch the overall strategy that has emerged:
\par
-- [1st step: Growth] A first crucial step, which was first
successfully taken by Helfgott~\cite{helfgott} for $G=\SL_2$ and
$\SL_3$~\cite{hsl3}, is to prove a growth theorem in the finite groups
$G(\Fp_p)$ for $p$ prime: there exists $\delta>0$, depending only on
$G$, such that for any generating subset $A\subset G(\Fp_p)$, we have
\begin{equation}\label{eq-growth}
|A\cdot A\cdot A| =|\{abc\,\mid\, a, b, c\in A\}|\gg
\min(|G(\Fp_p)|,|A|^{1+\delta}),
\end{equation}
where the implied constant depends only on $G$.
\par
Such a result, once known, implies that the diameter of the Cayley
graphs is $\ll (\log p)^C$ for some constant $C\geq 1$. This, by
itself, suffices -- by standard graph theory -- to obtain an explicit
upper-bound for the spectral radius $\rho_p$, but one which is weaker
than the desired uniform spectral gap (namely, of the type
$1-\rho_p\gg (\log p)^{-D}$ for some $D\geq 0$). Although this is
insufficient for applications to results like Theorem~\ref{th-bgs-ex},
it may be pointed out that this is enough for some others, including
rather surprising ones in arithmetic geometry~\cite{ehk}.
\par
After Helfgott's breakthrough, growth results were proved by Gill and
Helfgott~\cite{gill-helfgott} (for $\SL_m$, with a restriction on $A$)
and (independently and simultaneously) by
Breuillard-Green-Tao~\cite{bgt} and Pyber-Szab\'o~\cite[Th. 4]{psz} in
(more than) the necessary generality for our purpose.\footnote{\ This
  increased generality may be very useful for other applications,
  e.g., the Pyber-Szab\'o version is quite crucial in~\cite{ehk}.} An
intermediate paper of Hrushovski~\cite{hru} should be mentioned, since
it brought to light a somewhat old preprint of Larsen and
Pink~\cite{larsen-pink}, from which a useful general inequality
emerged (see, e.g,~\cite[Th. 4.1]{bgt}) concerning (roughly) the size
of the intersection of a ``non-growing'' set of $G(\Fp_p)$ and proper
algebraic subvarieties of $G$.
\par
-- [2nd step: Expansion for primes] As mentioned, the growth theorem
does not immediately imply that the Cayley graphs are
expanders. Bourgain and Gamburd~\cite{bg} were the first to prove this
property for $\SL_2(\Fp_p)$. Their method starts with an approach
going back to Sarnak and Xue~\cite{sarnak-xue}, which compares upper
and lower bounds for the number of loops in the Cayley graphs, based
at the identity, and of length $\ell\approx \log p$. As
in~\cite{sarnak-xue}, the lower bound comes from a spectral expansion
and the fact that the smallest degree of a non-trivial linear
representation of $\SL_2(\Fp_p)$ is ``large'' (namely, it is
$(p-1)/2$, as proved by Frobenius already). The upper-bound relies on
a new important and ingenious ingredient, now called ``flattening
lemma'' (\cite[Prop. 2]{bg}), which is used to show that large girth
of the Cayley graphs (a property which is fairly easy to prove) is
enough to ensure that, after $\gg \log p$ steps, the random walks on
the graphs are very close to uniformly distributed, and thus there
can't be too many loops of that length at the identity. In turn, the
proof of the flattening lemma turns out, ultimately, to be obtained
from Helfgott's growth result~(\ref{eq-growth}) in $\SL_2(\Fp_p)$.
(Why is that so? Very roughly, one may say that Bourgain and Gamburd
show that, if doubling the number of steps $\gg \log p$ of the random
walk does not lead to a great improvement of its uniformity, it must
be the case that its support must be to a large extent concentrated on
a set $A\subset \SL_2(\Fp_p)$ which does not grow, i.e., such
that~(\ref{eq-growth}) is false; according to Helfgott's theorem, this
means that $A$ is contained in a proper subgroup, but such a
possibility is in fact fairly easy to exclude, because one started
with a random walk using generators of $G(\Fp_p)$).
\par
After the proof of the general growth theorems, this second step was
extended to other groups (e.g., it is announced by Breuillard, Green
and Tao in~\cite{bgt}).
\par
-- [3rd step: Expansion for squarefree $d$] This step, which the
discussion above has shown to be absolutely fundamental for sieve
applications, was first done by Bourgain, Gamburd and Sarnak for
$\SL_2$ in~\cite{bgs}, using a rather sophisticated argument. However,
Varj\'u~\cite{varju} found a more streamlined proof, which can be
adapted to more general groups, in particular $\SL_m$, as soon as a
growth theorem for $G(\Fp_p)$ is known.\footnote{\ Note that Bourgain
  and Varj\'u~\cite{bv} also prove the expansion property for
  $\SL_m(\Zz/d\Zz)$ for all $d\geq 1$, not only those which are
  squarefree.} More general cases, including the statement we have
given, have been announced by Salehi Golsefidy and Varj\'u~\cite{sg-v}
(who give the most general, and in fact, best possible version, which
applies to any group such that the connected component of identity of
the Zariski-closure is \emph{perfect}) .

\subsection{Quantitative equidistribution: spectral and ergodic
  aspects}\label{sec-spectral}

We now consider a sieve in orbit, for a subgroup $\Lambda$ with
Zariski-closure $G/\Qq$, where we count using counting measures on
archimedean balls. The results are more fragmentary than in the
combinatorial case. Certainly, from the discussion above, we see that
the main issue is to extend the quantitative equidistribution
statement (Part (3) of Theorem~\ref{th-approx}) to this counting
method. However, Parts (1) and (2) still apply. Since
$$
\mu_X(\gamma\in\Lambda\,\mid\, \gamma\equiv\gamma_0\mods{d})=
\sum_{\stacksum{\|\gamma\|\leq X}{\gamma\equiv \gamma_0\mods{d}}}{1}
$$
for any $d\geq 1$ and $\gamma_0\in \Lambda_d$, and this is also
$$
\sum_{\stacksum{\|\tau\gamma_0\|\leq X}{\tau\in\Lambda(d)}}{1}
$$
where $\Lambda(d)=\ker(\Lambda\ra\Lambda_d)$ is the $d$-th
(generalized) congruence subgroup of $\Lambda$, one can see that this
amounts to issues of uniformity and effectivity in ``lattice-point
counting'' for the quotient $X_{\Lambda}=\Lambda\backslash G(\Rr)$ and
the congruence covers $X_{\Lambda}(d)=\Lambda(d)\backslash G(\Rr)$. 
\par
For $G(\Rr)=\SL_2(\Rr)$, $\Lambda\subset \SL_2(\Zz)$, and
$X_{\Lambda}$ of finite volume, a well-known result of Selberg (see,
e.g.,~\cite[Th. 15.11]{IK}), the original proof of which depends on
the spectral decomposition of the Laplace operator on $X_{\Lambda}$,
proves the local equidistribution with good error term depending
directly on the first non-zero eigenvalue $\lambda_1(d)$ for the
Laplace operator on the hyperbolic surface $\Lambda(d)\backslash
\mathbf{H}$. This indicates once more that spectral gaps -- of some
kind -- are crucial tools for the quantitative equidistribution. The
striking difference with the elementary argument leading
to~(\ref{eq-splitting}) should become clear: instead of counting
integers in a (large) interval, where the boundary contribution is
essentially negligible, we have hyperbolic lattice-point problems,
where the ``boundary'' may contribute a positive propertion of the
mass.
\par
It is now natural to distinguish two cases, depending on whether
$\Lambda$ is a lattice in the real points of its Zariski-closure $G$
(always assumed to be simple, connected and simply-connected), or
whether $\Lambda$ has infinite index in such lattices; in terms of
$X_{\Lambda}$, the dichotomy has very clear meaning: either
$X_{\Lambda}$ has finite or infinite volume, with respect to the
measure induced from a Haar measure on $G(\Rr)$. (Note that we still
require Theorem~\ref{th-approx} to be valid, which means that $G(\Rr)$
has no compact factor).
\par
(1) [Finite-volume case] Although it seems natural to apply methods of
harmonic analysis on $X_{\Lambda}(d)$, similar to Selberg's, there are
serious technical difficulties. This is especially true when
$X_{\Lambda}$ is not compact, since the full spectral decomposition of
$L^2(X_{\Lambda})$ depends then on the general theory of Eisenstein
series (see the paper of Duke, Rudnick and Sarnak~\cite{drs} for the
first results along these lines).
\par
However, starting with Eskin-McMullen~\cite{emm}, a number of methods
from ergodic theory have been found to lead to very general results on
lattice-point counting in this finite-volume case. For the purpose of
showing the required quantitative uniform equidistribution (as in
Theorem~\ref{th-approx}), one may mention first the results of
Maucourant~\cite{mau}; the most general ones have been extensively
developed by Gorodnik and Nevo~\cite{gn},~\cite{gn2} (see
also~\cite{n-s}). Without saying more (due to a lack of competence),
it should maybe only be said that the incarnation of the spectral gap
that occurs in this case is the exponent $p=p_{\Lambda}>2$ such that
matrix coefficients of unitary representations occurring in
$L^2_0(\Lambda(d)\backslash G(\Rr))$ are in $L^{p+\eps}$ for all
$\eps>0$. The existence of such a $p>2$ is known from the validity of
Property $(\tau)$. If $G(\Rr)$ has Property $(T)$, this constant
depends only on $G(\Rr)$, and explicit values are known (due to
Li~\cite{li} for classical groups and Oh in general~\cite{oh-p}); for
certain groups like $\SL_m$, $m\geq 3$ or $\Sp_{2g}$, $g\geq 2$, these
works give optimal values, as far as the general representation theory
of the group $G(\Rr)$ is concerned (the actual truth for congruence
subgroups lies within the realm of the Generalized Ramanujan
Conjectures, and is deeper; see~\cite{sarnak-ram} for more about these
aspects.)
\par
(2) [Infinite volume case\footnote{\ Bourgain, Gamburd and Sarnak say
  that this is the case of ``thin'' subgroup $\Lambda$, which is
  appealing terminology, but -- unfortunately -- clashes with the
  meaning of ``thin'' in Definition~\ref{def-thin} -- no Zariski-dense
  subgroup $\Lambda\subset \GL_m(\Zz)$ of $G$ is ``thin'' in
  $G(\Qq)$.}] The archimedean counting for these groups is the most
delicate among the cases currently considered. Indeed, the only
examples which have been handled in that case are subgroups of the
isometry groups of hyperbolic spaces, i.e., of orthogonal groups
$O(n,1)$, where the Lax-Phillips spectral approach to lattice-point
counting~\cite{lax-phillips} is available, at least when the Hausdorff
dimension of the limit set of the discrete subgroup $\Lambda\subset
\SO(n,1)(\Rr)$ is large enough.\footnote{\ Very recent work of
  Bourgain, Gamburd and Sarnak~\cite{bgs2} has started approaching the
  problem for subgroups of $\SL_2(\Rr)$ with limit sets of any
  positive dimension.} This is the case, for instance, for the
Apollonian group $\apol$ (which can be conjugated into a subgroup of
$O(3,1)(\Rr)$): the limit set has Hausdorff dimension $>1.30$, whereas
the Lax-Phillips lower-bound is $\delta>1$. 
\par
Again, due to a lack of competence, no more will be said about the
techniques involved, except to mention that the presence of a spectral
gap for the hyperbolic Laplace operator still plays a crucial role;
such gaps are established either by methods going back to Gamburd's
thesis~\cite{gamburd}, or by extending to infinite volume the
comparison theorems between the first non-zero eigenvalues for the
hyperbolic and combinatorial Laplace operators (due to Brooks and
Burger in the compact case, see, e.g.,~\cite[Ch. 6]{burger}), and
applying the corresponding case of expansion for Cayley graphs
(Theorem~\ref{th-expansion-groups}).  We will however state a few
results of Kontorovich and Oh~\cite{k-o} in Section~\ref{ssec-others},
and refer to Oh's ICM report~\cite{oh} for more on the methods
involved.

\subsection{Finiteness of saturation number}
\label{ssec-saturation}

We now show how to implement the sieve to prove
Theorem~\ref{th-bgs-ex}, using the weighted counting method (i.e.,
implicitly, random walks). It should be quite clear that the method is
very general.
\par
Let $Y=\Lambda$, $x_0$ and $f$ be as in the theorem, or indeed
Zariski-dense in a simple simply-connected group $G$ (instead of
$\SL_m$). For simplicity, we consider the sieve in $Y=\Lambda$ instead
of the orbit $\orbit(x_0)$; it is straightforward to deduce saturation
for the latter from this. We also assume that the irreducible
components of the hypersurface $\{f(\gamma x_0)=0\}$ in $G$ are
absolutely irreducible.\footnote{\. As noted in~\cite[p. 562]{bgs}, in
  the simply-connected case, the ring $\Qq[G]$ of functions on $G$ is
  factorial; the assumption is then that the irreducible factors of
  $f$ are still irreducible in $\bar{\Qq}[G]$.}  Fix a symmetric set
of generators $S$ with $1\in S$. By our previous arguments
(Theorem~\ref{th-approx} and Theorem~\ref{th-expansion-groups}), the
basic requirements of sieve are met.
\par
We proceed to study the sifted set~(\ref{eq-sifted}) for the set of
primes $\mathcal{P}$ consisting of those $p$ not in the finite
``exceptional'' set $\sop$ given by Theorem~\ref{th-approx}, and
$$
\Omega_p=\{\gamma\in Y_p=G(\Fp_p)\,\mid\, f(\gamma\cdot
(x_0\mods{p}))=0\in \Zz/p\Zz\}.
$$
\par
Indeed, $\sifted_z(\Lambda;\Omega)$ is the set of $\gamma\in\Lambda$
for which $f(\gamma\cdot x_0)$ has no prime factor $<z$ (outside
$\sop$). After maybe enlarging the set $\sop$ (remaining finite),
standard Lang-Weil estimates show that
$$
|\Omega_p|=\kappa p^{\dim(G)-1}+O(p^{\dim(G)-3/2})
$$
where $\kappa$ is the number of absolutely irreducible components of
the hypersurface in $G$ defined by $\{f(\gamma x_0)=0\}$, and the
implied constant is absolute. Since
$$
|G(\Fp_p)|=p^{\dim(G)}+O(p^{\dim(G)-1/2}),
$$
(which can be checked very elementarily for many groups) it follows
that the density of $\Omega_p$ satisfies
$$
\nu_p(\Omega_p)=\frac{\kappa}{p}+O(p^{-3/2})
$$
for all $p\notin\sop$. This verifies~(\ref{eq-dim-kappa}): the sieve
in orbit has ``dimension'' $\kappa$ in the standard sieve
terminology. We see from~(\ref{eq-delta})\footnote{\
  Remark~(\ref{rm-finitude}) is applicable here to
  check~(\ref{eq-fin-unite}).} that~(\ref{eq-easy-sieve}) holds with
$$
D_k=\beta^k
$$
for some $\beta>1$ (indeed $\beta$ can be any real number
$<\rho^{-1/(1+3\dim(G)/2)}$ , where $\rho<1$ is the expansion constant
for our Cayley graphs, as in Definition~\ref{def-expand}).
\par
The conclusion is that there are many $\gamma\in\Lambda$ where
$f(\gamma\cdot x_0)$ is not divisible by primes $<z=\beta^{k/s}$,
indeed the $\mu_k$-measure of this set, say $\sifted_k$, is
$$
\mu_k(\sifted_k)\gg \frac{1}{(\log z)^{\kappa}}\asymp 
\frac{1}{k^{\kappa}},
$$
for $k$ large enough. To prove from this that the saturation number is
finite, we need two more easy ingredients:
\par
(1) If $\gamma\in\sifted_k$, then the integer $n=f(\gamma\cdot x_0)$
has a \emph{bounded} number of prime factors. (Except if $n=0$, which
only happens with much smaller probability, see
Remark~\ref{rm-finitude}). Indeed, we have
$$
n=f(s_1\cdots s_n x_0)
$$
for some $s_i\in S$. Since the function $f$ has polynomial growth, we
see immediately that there exists a constant $\lambda\geq 1$ such that
\begin{equation}\label{eq-apriori}
f(\gamma\cdot x_0)\ll \lambda^k
\end{equation}
for all $\gamma\in \sifted_k$. An integer of this size, with no prime
factor $<\beta^{k/s}$, must necessarily satisfy
$$
\Omega(f(\gamma\cdot x_0))\leq r=\frac{s\log\lambda}{\log\beta}.
$$

\begin{rema}
  If the Cayley graphs satisfy a weaker property than expansion, one
  can still do a certain amount of sieving. However, the level of
  distribution being weaker, one obtains only points of the orbit with
  fewer prime factors than typically expected for integers of that
  size (see the Appendix for the meaning of this).
\end{rema}
\par
(2) We must check that the lower bound is incompatible with the set of
$\gamma$ with $\Omega(f(\gamma\cdot x))\leq r$ being too small, i.e.,
thin or simply not Zariski-dense. For the latter this is quite easy:
indeed, any subset $W$ of $\Lambda$ contained in a proper hypersurface
$\{g=0\}$ of $G$ satisfies the much slower growth
$$
\mu_k(W)\ll \delta^{-k}
$$
for some $\delta>1$, as one can see simply by selecting a suitable
prime $p$ for which $\{g=0\}$ is a hypersurface modulo $p$ and
bounding
$$
\mu_k(W)\leq \mu_k(\gamma\,\mid\, g(\gamma)=0\mods{p})
$$
using local equidistribution modulo $p$ and the Lang-Weil estimates. 
\par
In order to show a similar result for thin sets, however, one must
apply the large sieve instead, as discussed briefly in
Section~\ref{sec-ls}.

\subsection{Other results for the sieve in orbits}
\label{ssec-others}

We collect here a few results which have been proved in the setting of
the sieve in orbits.

\begin{exem}
  We start with results concerning the Apollonian group and the
  associated circle packings.   
\begin{itemize}
\item Fuchs~\cite{fuchs} has studied very carefully the reductions
  modulo integers of the Apollonian group.
\item Based on this study, a delicate conjecture predicts a
  local-global principle for the presence of integers among the
  curvature set $\curv(\uple{c})$; Bourgain and Fuchs~\cite{b-f} have
  at least shown that the number of integers $\leq T$ arising as
  curvatures (without multiplicity) is $\gg T$.
\item Kontorovich and Oh have applied spectral-ergodic counting
  methods in infinite volume to deduce, first, asymptotic formulas for
  the number of curvatures $\leq T$,\footnote{\ Counted with
    multiplicity; the latter, on average, is quite large: about
    $T^{\delta-1}$ where $\delta>1.3$ is the dimension of the limit
    set.} and then -- by means of sieve -- have obtained upper and
  lower bounds for the number of prime curvatures, or the number of
  pairs of prime curvatures of two tangent circles in the packing
  (such as $11$ and $23$ in Figure~1). Note that here, counting in the
  orbit means that the Lax-Phillips theory does not apply, and thus
  new ideas are needed. They also did a similar analysis for orbits of
  infinite-index subgroups $\Lambda$ of $\SO(2,1)(\Zz)$ acting on the
  cone of Pythagorean triples (see Remark~\ref{rm-pythagor}),
  see~\cite{k-o-2}; remarkably, using Gamburd's explicit spectral
  gap~\cite{gamburd}, they obtain for instance -- for sufficiently
  large limit sets, but possiby infinite index -- the expected
  proportion of triangles with hypothenuse having $\leq 14$ prime
  factors.
\end{itemize}
\end{exem}

\begin{exem}\label{ex-vmn}
  The integral points $V_{m,n}$ of the $\SL_m$-homogeneous spaces
$$
\mathcal{V}_{m,n}=\{\gamma\in \GL_m\,\mid\, \det(\gamma)=n\}
$$
have been studied in great detail by Nevo and Sarnak~\cite{n-s}, in
the setting of archimedean balls, using methods based on mixing and
ergodic theory.  They show, for instance, that if $f\in
\Qq[\mathcal{V}_{m,n}]$ is integral valued on $V_{m,n}$, absolutely
irreducible and has no congruence obstruction to being prime (i.e.,
for any prime $p$, there exists $\gamma\in V_{m,n}$ with $p\nmid
f(\gamma)$), then the saturation number of $V_{m,n}$ is
$$
\leq 1+18m_e^3\deg(f),
$$
where $m_e$ is the smallest even integer $\geq m-1$, in fact that
\begin{equation}\label{eq-satur-arch}
|\{\gamma\in V_{m,n}\,\mid\, \|\gamma\|\leq T\text{ and }
\Omega(f(\gamma))\leq r\}|\gg 
\frac{
|\{\gamma\in V_{m,n}\,\mid\, \|\gamma\|\leq T\}|}{(\log T)},
\end{equation}
for $r>18m_e^3\deg(f)$.
\end{exem}

\begin{exem}
  As explained in~\cite{n-s}, bounds like~(\ref{eq-satur-arch}) do not
  transfer trivially to non-principal homogeneous spaces (i.e., orbits
  of an arithmetic group with non-trivial stabilizer), although this
  is no problem when the mere finiteness of a saturation number is
  expected. Gorodnik and Nevo~\cite{gn, gn2} have obtained results
  which extend such results to many cases. Their results apply, for
  example, to the orbits
$$
\orbit(g_0)=
\{g\in M_m(\Zz)\,\mid\, g=\T{\gamma}g_0\gamma\text{ for some }
\gamma\in \SL_m(\Zz)\}
$$
for a fixed non-degenerate symmetric integral matrix $g_0$, if $m\geq
3$. Thus, for suitable functions $f$, $\kappa$ and (explicit) $r$,
they prove a lower bound
$$
|\{g\in \orbit(g_0)\,\mid\, \|g\|\leq T\text{ and }
\Omega(f(g))\leq r\}|\gg 
\frac{
|\{g\in \orbit(g_0)\,\mid\, \|g\|\leq T\}|}{(\log T)^{\kappa}}
$$
(here the stabilizer is an orthogonal group).
\end{exem}

\section{Related sieve problems and results}

We present here other developments of sieve in expansion, as well as
some analogues over finite fields.

\subsection{Geometric examples}\label{sec-geom}

In the spirit of Section~\ref{ssec-geom}, there are a number of
geometric situations where one naturally wonders about ``genericity''
properties of elements in interesting discrete groups not given as
subgroups of some $\GL_m(\Zz)$. Sometimes, using arithmetic quotients,
and their reductions modulo primes, is enough to attack very
interesting problems, as we described already for the homology of
Dunfield-Thurston $3$-manifolds.
\par
For these, the discrete group involved is the mapping class group
$\Gamma_g$ of a surface $\Sigma_g$ of genus $g$. It is finitely
generated, and because rather little is known about the precise
structure of combinatorial balls, it is natural (as done in~\cite{dt})
to use a weighted combinatorial counting to apply sieve in that case,
or in other words, to use a random walk on $\Gamma_g$ based on a
symmetric generating set $S$, with $1\in S$ for simplicity. 
\par
Since~(\ref{eq-homology}) and~(\ref{eq-homologyp}) only depend on the
image of a mapping class $\phi$ in $\Sp_{2g}(\Zz)$ or
$\Sp_{2g}(\Fp_p)$, one is -- in effect -- doing a random walk (though
not always with uniformly probable steps) on the discrete group
$\Sp_{2g}(\Zz)$. Since, for $g\geq 2$ (which is most interesting) this
group has Property (T), the basic requirements of sieve hold.
\par
In addition, it is not difficult to compute the size of
$$
\Omega_p=\{\gamma\in \Sp_{2g}(\Fp_p)\,\mid\, \langle J_p,\gamma J_p\rangle
\not=\Fp_{p}^{2g}\},
$$
which is of size $p^{-1}+O(p^{-2})$ for $p\geq 2$ (for fixed $g$;
intuitively a determinant must be zero for this to hold, and this
happens with probability roughly $1/p$). Thus the homology of
Dunfield-Thurston $3$-manifolds can be handled with a sieve of
dimension $1$. 
\par
If $\phi_k$ is the $k$-th step of a random walk on $\Gamma_g$ (with
respect to a generating set), using notation from
Section~\ref{ssec-local-global}, the set $Y^0$ corresponds to those
manifolds with $H_1(M_{\phi_k},\Zz)$ which is infinite. As in
Remark~\ref{rm-finitude}, this event has probability to $0$ (proved
in~\cite{dt}) exponentially fast as $k\ra +\infty$ (\cite[Pr. 7.19
(1)]{cup}).
\par
One can then also deduce that
$$
\proba(H_1(M_{\phi_k},\Zz)\text{ has no $p$-part for } p<z=\beta^{k})
\asymp \frac{1}{k},
$$
for some $\beta=\beta(g)>1$. Using the
description~(\ref{eq-homology}), we see also that if
$H_1(M_{\phi_k},\Zz)$ is finite, its order can not be too large, more
precisely there exists $\lambda\geq 1$ such that the product
$\Delta_k$ of those $p$ with $H_1(M_{\phi_k},\Fp_p)\not=0$ satisfies
$$
p\leq \lambda^k
$$
(because $\Delta_k$ divides a non-zero determinant of such size). So
by comparison, we deduce that there exists $r$ (depending on $g$ and
the generators $S$) such that
$$
\proba(H_1(M_{\phi_k},\Zz)\text{ is finite and has order divisible by
  $\leq r$ primes})\gg \frac{1}{k}.
$$
\par
In another direction, an application of the large sieve shows that,
with probability going to $1$, $|H_1(M_{\phi_k},\Zz)|$ is divisible by
``many'' primes $<z$. This means $|H_1(M_{\phi_k},\Zz)|$ is typically
finite, but very large (see~\cite[Pr. 7.19 (2)]{cup}).
\par
Other examples of groups where sieve can be applied are given by
automorphisms of free groups (of rank $m\geq 2$, where the action on
the abelianization gives a quotient $\SL_m(\Zz)$).  Thus, sieve
methods give another illustration of the many analogies between these
discrete groups (others are surveyed in the recent talk~\cite{paulin}
of F. Paulin in this seminar); see also the related works of
Rivin~\cite{rivin} and Maher~\cite{maher}.


\subsection{Large sieve problems}\label{sec-ls}

We have briefly mentioned the large sieve already, and we will now add
a few words (see~\cite{cup} for much more on this topic). In the
context of Section~\ref{ssec-local-global}, and starting with some
work of Linnik, many sieve situations have appeared where the
condition sets $\Omega_p$ satisfy
\begin{equation}\label{eq-large}
\nu_p(\Omega_p)\geq \delta>0
\end{equation}
for some $\delta>0$ and all primes $p$; because, for the classical
case, this amounts to excluding many residue classes, it is customary
to speak of a \emph{large sieve} situation.
\par
The large sieve method, under the assumptions of Basic Requirements
(local independent equidistribution and its quantitative version)
leads roughly to two types of statements:\footnote{\ Where it is not
  necessary to assume, a priori, that~(\ref{eq-large}) holds.}
\par
(1) An upper-bound for $\mu_n(\sifted_z(Y;\Omega))$ of the type
$$
\mu_n(\sifted_z(Y;\Omega))\ll \mu_n(Y)H^{-1},\quad\quad
H=\sum_{d<z}{\mu(d)^2\prod_{p\mid
    d}{\frac{\nu_p(\Omega_p)}{1-\nu_p(\Omega_p)}}}
$$
(for $z$ of size similar to the level of distribution;
see~\cite[Prop. 2.3, Cor. 2.13]{cup}). In the sieve in orbits, this
can be used to show that
$$
\mu_n(W)\ll \delta^{-k}
$$
for some $\delta>1$, where $W\subset G(\Qq)\cap \GL_m(\Zz)$ is a thin
set, using the fact (see~\cite[Th. 3.6.2]{serre-galois}) that the
complement $\Omega_p$ of $W\mods{p}$ satisfies a large-sieve
condition:
$$
|\Omega_p|\gg 1
$$
for $p$ large enough. This extends the finiteness of saturation
numbers to thin sets.
\par
(2) An upper-bound for the mean-square of
$$
\Bigl(\sum_{\stacksum{p<z}{y\mods{p}\in \Omega_p}}{1}-
\sum_{p<z}{\nu_p(\Omega_p)}\Bigr)
$$
with respect to $\mu_n$ (see~\cite[Prop. 2.15]{cup}). This leads to
the fact that the number of $p<z$ such that $y\mods{p}$ is in
$\Omega_p$ is close to the expected value
$$
\sum_{p<z}{\nu_p(\Omega_p)}
$$
with high probability.
\par
This is used for instance to show that the homology of the
$3$-manifolds has typically a very large torsion part (growing faster
than any polynomial, as $k\ra +\infty$).
\par
Another application of the large sieve, in settings related to
discrete groups, concernes the question of trying to detect the
``typical'' Galois group of the splitting field of the characteristic
polynomial of an element $x$ in a subgroup $\Lambda\subset
\GL_m(\Zz)$. The (quite classical) idea is to use Frobenius
automorphisms at primes to produce conjugacy classes in the Galois
group. If, for instance, $\Lambda$ is Zariski-dense in $\SL_m$, the
factorization pattern of the characteristic polynomial modulo $p$
gives a conjugacy class $c$ in the symmetric group $\mathfrak{S}_m$,
which is the maximal possible Galois group for the characteristic
polynomial. Since it is not too difficult to show that
$$
|\{g\in \SL_m(\Fp_p)\,\mid\,
\text{the conjugacy class associated to $g$ is $c$}\}|\sim
\frac{|c|}{|\SL_m(\Fp_p)|}
$$
for fixed $m$, conjugacy class $c$ and $p\ra +\infty$, this is a
condition like~(\ref{eq-large}) when $\nu_p$ is the uniform
probability measure on $\SL_m(\Fp_p)$. One can then prove that the
Galois group is as large as possible, with probability going to $1$
(see~\cite{rivin},~\cite[Th. 7.12]{cup} and, for a very general
statement, the recent work of the author with F. Jouve and
D. Zywina~\cite{jkz}, where the typical Galois group is essentially
the Weyl group of the Zariski closure of $\Lambda$.)
\par
Finally, very recently, Lubotzky and Meiri~\cite{lm} have used the
large sieve (and the expansion result of Salehi Golsefidy and
Varj\'u~\cite{sg-v}) in order to prove that, if $\Gamma$ is a
finitely-generated subgroup of $\GL_m(\Cc)$ which is not virtually
solvable, one has
$$
\proba(X_k\text{ is of the form $\gamma^m$ for some $\gamma\in\Gamma$
  and $m\geq 2$})
\ll \exp(-\beta k)
$$
for every left-invariant random walk $(X_k)$ on $\Gamma$ defined using
a symmetric generated set $S$ of $\Gamma$ (with $1\in\Gamma$), where
$\beta>0$ depends on $S$. This result does not have any obvious
``classical analogue'', and is a strong form of a converse of a result
of Mal'cev. Moreover, the proof involves many subtle group-theoretic
ingredients in addition to the sieve, and hence this theorem seems to
be an excellent illustration of the potential usefulness of sieve
ideas as a new tool in the study of discrete groups.

\subsection{Sieve for Frobenius over finite fields}

There are a number of interesting analogies between the type of sieve
problems in Section~\ref{sec-orbits} and problems of arithmetic
geometry over finite fields which concern the properties of the action
(typically, characteristic polynomials) of Frobenius elements
associated to families of algebraic varieties over finite fields. In
this context, instead of expansion properties, one uses the Riemann
Hypothesis over finite fields (and uniform estimates for Betti
numbers) to prove the required equidistribution properties
(quantitative uniform versions of the Chebotarev density theorem). We
refer to~\cite[\S 8, App. A]{cup} for precise descriptions and sample
problems, especially in large-sieve situations (which were already
implicit in work of Chavdarov~\cite{chavdarov}), and only mention that
a prototypical question is the following: given $f\in\Fp_p[X]$ of
degree $2g$ and without repeated roots, and the family of
hyperelliptic curves given by
$$
C_t\,:\, y^2=f(x)(x-t)
$$
where $t$ is the parameter, how many $t\in \Fp_{p^{\nu}}$ are there
such that $|C_t(\Fp_{p^{\nu}})|$ is prime, or almost prime?
\par
One may also consider the case of a fixed algebraic variety over a
number field, and the variation with $p$ of its reductions modulo
primes. The principles of the ``sieve for Frobenius'' (now, in some
sense, in horizontal context) are still applicable, though they suffer
from the lack of Riemann Hypothesis (or even strong enough versions of
the Bombieri-Vinogradov Theorem), and the unconditional results are
therefore  fairly weak (see~\cite{chavdarov} and~\cite{zywina}).

\section{Remarks, problems and conjectures}\label{sec-conclusion}

We conclude by describing some open interesting problems and other
related works.
\par
(1) [Conjugacy classes] Let $\Lambda\subset \SL_m(\Zz)$ be a
Zariski-dense subgroup, of infinite index. The theory and results
described previously give information -- theorems or conjectures --
concerning the distribution of the elements of $\Lambda$ and, in some
sense, their density among the elements of $\SL_m(\Zz)$. Now one may
ask: what about the set of conjugacy classes of $\Lambda$?  This seems
like a very natural and interesting question, and even in the case of
$m=2$ it does not seem (to the author's knowledge) that much is
known. 
\par
(2) [Strong equidistribution for word-length metric] It would be of
great interest to obtain a version of Part (3) of
Theorem~\ref{th-approx} for the (unweighted) word-length counting
method when $\Lambda$ is a fairly general group with exponential
growth (in particular, when it is not free). 
\par
(3) [Explicit bounds] We have concentrated on general results, which
in some sense are quite basic from the point of view of applying
sieve. It is natural that now much effort goes into improving the
results, and in particular in obtaining explicit bounds for saturation
numbers,\footnote{\ Where ``explicit'' means having a concrete number,
  be it $10$, $100$ or $1000$ for a concrete case like, for instance,
  the group $L$ and the polynomial $f(\gamma)=\text{product of the
    coordinates}$.} or explicit quantitative lower-bounds. We have
mentioned examples like those of Nevo-Sarnak or Gorodnik-Nevo for
lattices and archimedean balls. These, as well as the argument in
Section~\ref{ssec-saturation}, indicate clearly that a first
inevitable step is to prove an explicit version of spectral gap. In
the ergodic setting, this comes ultimately from spectral theory, and
the gap is quite explicit (this goes back to Selberg's famous $3/16$
theorem).
It would be extremely interesting to have, for instance, a version of
Theorem~\ref{th-expansion-groups} (even, to begin with, for $\SL_2$)
in which the expansion constant for the Cayley graph is a known
function of, say, the coordinates of the matrices in the generating
set $S$.  As pointed out by E. Breuillard, the issue is not the
effectiveness of the methods (for instance, there is no issue
comparable to the Landau-Siegel for zeros Dirichlet $L$-functions):
the methods and results that lead to this theorem are effective in
principle (but one must be careful when general groups are involved
and ``effective'' algebraic geometry is needed).
\par
(4) [Refinements] Once -- or when -- an explicit spectral gap is
known, one can envision the application of more refined versions of
sieve; this has been done, e.g., by Liu and Sarnak~\cite{liu-s} for
sieving integral points on quadrics in three variables, where a
sophisticated weighted sieve is brought to bear.
\par
Along these lines, it would be extremely interesting also to find
examples where Iwaniec's Bilinear Form of the remainder term for the
linear sieve (see~\cite[\S 12.7]{FI}) was exploited. Similarly, it
would be remarkable to have applications where the level of
distribution is obtained by a non-trivial average estimate of the
remainders $r_d$, instead of summing individual esimates (this being
the heart of the Bombieri-Vinogradov theorem).
\par
(5) [Primes?] In many cases, when there are no congruence
obstructions, one expects that the saturation number be $1$, i.e.,
that many elements of an orbit have $f(x)$ be prime\footnote{\ More
  precisely, prime or opposite of a prime; for rather fundamental
  reason, explained in~\cite[\S 2.3]{bgs}, one cannot hope to
  distinguish between these two.}. For instance, Bourgain, Gamburd and
Sarnak propose~\cite[Conjecture 1.4]{bgs} a fairly general conjecture
concerning the value of the saturation number. The paper of Fuchs and
Sanden~\cite[Conj. 1.2, 1.3]{f-s} gives two very precise quantitative
conjectures concerning prime curvatures of Apollonian circle packings,
which are quite delicate (and shows that making quantitative
conjectures is rather subtle in such settings). Some results with
primes are known, but the methods used are more directly comparable
with those of Vinogradov and the circle method than with sieve: Nevo
and Sarnak~\cite[Th. 1.4]{n-s} find a Zariski-dense subset of
$V_{m,n}$ (see~(\ref{ex-vmn})) where all coordinates of the matrix are
primes (up to sign), under the necessary condition that $n\equiv
0\mods{2^{m-1}}$, and Bourgain and Kontorovich~\cite{bk} show that
(for instance) the set of all integers arising as (absolute value of)
the bottom-right corner of an element in a thin subgroup of
$\SL_2(\Zz)$ with sufficiently large limit set contains all positive
integers $\leq N$ with $\ll N^{1-\delta}$ exceptions -- in particular,
infinitely many primes -- for $N$ large enough.
\par
One may then also ask: ``what is the strength of such statements, if
valid''?  What do they mean about prime numbers?  The only clue in
that direction -- to the author's knowledge -- is the following
indirect fact: Friedlander and Iwaniec have shown (see~\cite[\S
14.7]{FI}) that one can prove that there is the expected proportion of
matrices
$$
g=\begin{pmatrix}a&b\\c&d\end{pmatrix}
\in \SL_2(\Zz)
$$ 
with
$$
a^2+b^2+c^2+d^2=p\text{ prime},\ p\leq X,
$$
\emph{provided} a suitable form of the Elliott-Halberstam conjecture
holds (level of distribution $Q=X^{1/2+\delta}$, for some small
$\delta>0$, for primes $\leq X$). This is a type of sieve in orbits,
obviously. The assumption here is widely believed to be true, but
seems entirely out of reach (e.g., it is much stronger than the
Generalized Riemann Hypothesis). It is now known (by the work of
Goldston, Pintz and Y\i ld\i r\i m) that this would also imply the
existence of infinitely gaps of bounded size between consecutive
primes (see, e.g.,~\cite[Th. 7.17]{FI}).

\section*{Appendix: What to expect from integers}\label{sec-final} 

We recall here, very briefly and for completeness, the most basic
estimates concerning multiplicatives properties of integers. These
serve as comparison points for statements on the distribution of prime
factors of elements of any set of integers. Of course, all these facts
are known in much stronger form than what we state.

\begin{itemize}
\item The number of primes $p\leq X$ is asymptotic to $X/(\log X)$
  (the Prime Number Theorem).
\item More generally, for $k\geq 1$ (fixed), the number of integers
  $n\leq X$ which are product of $k$ (or $\leq k$) prime factors, is
  asymptotic to
$$
\frac{1}{(k-1)!}\frac{X(\log\log X)^{k-1}}{(\log X)}.
$$
\item On the other hand, for $k\geq 1$ (fixed), the number of integers
  $n\leq X$ which have \emph{no prime factor} $p\leq X^{1/k}$ is of
  order $\asymp X/(\log X)$. Note that this set is a subset of the
  previous one (when $\leq k$ prime factors are considered), but the
  restriction on the size of prime factors is more stringent than the
  restriction on their numbers, and the order of magnitude becomes
  insensitive to the number of prime factors $k$.
\item The typical number of prime divisors of an integer $n\leq X$ is
  $\log\log X$; indeed, we have the Hardy-Ramanujan variance bound
$$
\sum_{n\leq X}{\Bigl(\Omega(n)-\log\log X\Bigr)^2}
\ll X\log\log X, 
$$
so that, e.g., there are
$$
\ll \frac{X}{\log\log X}
$$
integers $\leq X$ with $|\Omega(n)-\log\log X|\geq (\log\log X)/2$. 
\end{itemize}

\end{document}